\documentclass[11pt,a4paper]{article}
\usepackage{amssymb,fullpage,amsmath,mathrsfs}
\usepackage{pgfplots}
\usepackage{amsfonts,amsbsy,graphicx}
\usepackage{authblk,lineno}
\usepackage{bm}
\usepackage{bbm}
\usepackage{mathtools}

\newcommand{\eref}[1]{(\ref{eqn:#1})}
\newcommand{\elab}[1]{\label{eqn:#1}}
\newcommand{\fref}[1]{\ref{fig:#1}}
\newcommand{\flab}[1]{\label{fig:#1}}
\newcommand{\sref}[1]{\ref{sec:#1}}
\newcommand{\slab}[1]{\label{sec:#1}}

\def\beq{\begin{equation}}
\def\eeq{\end{equation}}

\def\eps{\epsilon}

\def\XXint#1#2#3{{\setbox0=\hbox{$#1{#2#3}{\int}$}
\vcenter{\hbox{$#2#3$}}\kern-.5\wd0}}

\makeatletter
\newenvironment{taggedsubequations}[1]
 {%
  \addtocounter{equation}{-1}%
  \begin{subequations}%
  \def\@currentlabel{#1}%
   \renewcommand{\theequation}{#1\alph{equation}}%
 }
 {\end{subequations}}
\makeatother 

\begin{document}

\title{KPP fronts in shear flows with cut-off reaction rates}

\author{D. J. Needham and A. Tzella\thanks{Address for correspondence:   Prof. D. J. Needham and Dr A. Tzella, School of Mathematics, University of Birmingham; email:  a.tzella@bham.ac.uk}} 
\affil{School of Mathematics, University of Birmingham, Birmingham, B15 2TT, UK.}
\maketitle

\begin{abstract}
We consider the effect of a shear flow which has, without loss of generality, a zero mean flow rate, on a Kolmogorov–Petrovskii–Piscounov (KPP) type model in the presence of a discontinuous cut-off   at concentration $u = u_c$.
In the long-time limit, a permanent form travelling wave solution is established which, for fixed $u_c>0$, is unique.  Its structure and speed of propagation depends on $A$ (the strength of   the flow relative to the  propagation speed in the absence of advection) and $B$ (the square of the front thickness relative to the channel width). 
We use matched asymptotic  expansions to approximate the propagation 
  speed  in the three natural cases $A\to \infty$, $A\to 0$ and $A=O(1)$, with particular associated orderings on $B$, whilst $u_c\in(0,1)$ remains fixed. 
In all the cases that we consider, the shear flow enhances the speed of propagation in a  manner that is similar to the case without cut-off ($u_c=0$). 
We illustrate the theory by evaluating expressions (either directly or through numerical integration) for  the particular cases of the plane Couette and Poiseuille flows.

\end{abstract}

\section{Introduction}
We investigate travelling fronts for cut-off  reaction--diffusion equations evolving in an  infinite   channel in the presence of a shear flow.  The model equation that we consider describes the spatio-temporal evolution of the scalar function $u$ denoting the concentration field  of a dissolved  species, and
 takes the non-dimensional generic form  
 \begin{taggedsubequations}{IBVP}\label{IBVP}
\beq\elab{KPP}
u_t +  A\, \alpha(y) u_x=  u_{xx}+  B \, u_{yy}+
 f_c(u),~~
 (x,y)\in \mathbb{R}\times(0,1),~t\in\mathbb{R^+}.  
  \eeq
  {Here, the dimensionless   spatial coordinate    $y$ has been  obtained by scaling with  the channel width $a$ whilst 
  the dimensionless   spatial coordinate    $x$ has been obtained by scaling with the diffusive length scale $\sqrt{\kappa\tau}$ where $\kappa$ is the molecular diffusivity and $\tau$ is the 
   reaction timescale. 
  The latter has also been used to introduce the  
  dimensionless time $t$. 
  }
	  The cut-off reaction function
	    $f_c: \mathbb{R} \to \mathbb{R}$ is 
		given by
		 \beq \elab{BDreaction}
		 f_c(u)=
		 \begin{cases}
		   			f(u), & \text{$u >u_c$}\\
		             0, & \text{$u\leq u_c$}
		 		 \end{cases}
		\eeq
		where $f(u)$ is a normalised KPP-type reaction function, named after the
		classical work by  {Fisher \cite{Fisher1937}} and Kolmogorov, Petrovskii, and Piskunov \cite{Kolmogorov_etal1937}, 
which satisfies  
	  \beq
	  \begin{split}
			\elab{KPPreaction}
	  	  	 f(0) = f(1)= 0,\quad f'(0)  = 1, \quad & f'(1)  <0 \quad\text{and} 
			\quad f \in C^1(\mathbb{R}) \\
	  	  	 0<f(u)\leq u   \quad  \forall    u \in (0,1),\quad
	  	 	\quad & f(u) < 0 \quad  \forall    u \in (1,\infty).
			\end{split}
	  	 	\eeq
      Thus the reaction is effectively deactivated at points where the concentration $u$ lies at or below a threshold value $u_c \in (0, 1)$.
	  {For convenience throughout, we introduce the shorthand notation $f_c\coloneqq f_c(u_c^+)$, the limiting value from above of $f_c(u)$ at $u=u_c$. }
	  	  A prototypical example of such a KPP reaction function is the Fisher reaction function \cite{Fisher1937} given by  $f(u)=u(1-u)$.
		  The flow function $\alpha:[0,1]\to \mathbb{R}$ satisfies $\alpha\in C^1([0,1])$. It corresponds to 
  the $x$-component of a unidirectional, zero-mean, steady incompressible shear flow with fluid velocity field
		 \beq\elab{eqn1.4}
		  q(x, y) =   (\alpha(y),0) ~~\text{for}~~(x,y)\in\mathbb{R}\times[0,1], \quad\text{with}\quad \int_0^1 \alpha(y)dy=0. 
		  \eeq
	Equation \eref{KPP} involves two non-dimensional parameters 
 \beq
 A=  \frac{V}{\sqrt{\kappa/\tau}}    
 \quad
 \text{and}
 \quad
B=\frac{\kappa \tau}{a^2}, \nonumber
\eeq
 where   
 $V$ is  the characteristic speed  of the (dimensional) flow.
 They  respectively  
 measure the strength of advective fluid velocity relative to the    front speed (in the absence of advection),  and
  the square of the  front thickness (in the absence of advection) to the square of the channel width.
 We focus on two-dimensional straight infinite channel domains with Neumann boundary conditions at  the domain walls, 
 \beq\elab{.noflux}
 u_y(x,0)=u_y(x,1)= 0 \quad\text{for}~~x\in \mathbb{R}
 \eeq
 so that there is no flux of $u$  
 across the wall of the domain. 
These are supplemented by the  initial conditions and far field boundary conditions, namely,
    \beq\elab{IBCb} %
     u(x,y,0)=H(-x) ~~\text{for}~~ (x,y)\in \mathbb{R}\times[0,1],\quad  u(x,y,t) \to  
 	 \begin{cases}
 		 1, & \text{as $x \to - \infty$}\\
 		             0, & \text{as $x \to \infty$}
 	\end{cases}				 
  \eeq
   with these limits being uniform for  $y\in[0,1]$ and $t\in[0,T]$, for any $T>0$. Here $H:\mathbb{R}\to\mathbb{R}$ is the usual Heavyside function.  
   {In what follows we will refer to the initial boundary value problem 
   specified by  the above equations 
   as (IBVP)}.
\end{taggedsubequations}
   \subsection{Background and related works}

The (IBVP) was first considered  by
Brunet and Derrida \cite{BrunetDerrida1997} in the absence of a background flow {that is, when $A=0$}.  
They   proposed it as  a  model  to provide insight into    
 the spreading of discrete systems of interacting  chemical and biological particles
 in homogeneous environments. 
They conjectured that discreteness in concentration values can be represented by an effective cut-off concentration $u_c$, where $u_c$ may be viewed as the effective mass of a single particle. 
The idea is that for   $u < u_c$, diffusion dominates over growth. 
It is now recognised, following previous work by Tisbury and the authors \cite{Tisbury_etal2020b}, that the solution to (IBVP) in the absence of a {background}  flow  evolves at large times to a  permanent form traveling wave (PTW) solution $u(x,y,t)=U_T(x-vt)$ 
 and propagation speed $v=v^*(u_c)$ in the sense that there exist functions $s:\mathbb{R}^+\to \mathbb{R}$ and $U_T:\mathbb{R}\to \mathbb{R}$ such that
 \beq\elab{lim}
 s(t)t^{-1}\to v^*(u_c) \quad\text{and}\quad u(z+s(t),y, t)\to U_T(z),
 \eeq
as $t\to \infty$ uniformly for $z\in\mathbb{R}$, with 
 $(\dot s(t)-v^*(u_c))$ being exponentially small in $t$ as $t\to\infty$ and a similar, spatially uniform, rate of convergence in the second limit in \eref{lim}. Tisbury and the authors \cite{Tisbury_etal2020a} showed that the PTW solution  is, for each $u_c\in (0,1)$ and KPP function $f(u)$,  unique  (up to translation), monotone decreasing and positive, with $\lim_{z\to-\infty}U_T(z)=1$ and $\lim_{z\to\infty}U_T(z)=0$.  
Its propagation speed $v=v^*(u_c)$ is, for fixed $u_c\in(0,1)$, uniquely determined by the existence of a heteroclinic connection in the $(U,U')$ phase plane which connects the equilibrium point $(1,0)$, as $z\to-\infty$, to the equilibrium point $(0,0)$ as $z\to\infty$ (the translational invariance is then fixed by  requiring that  $U_T(0)=u_c$).  
 An explicit expression for  $v^*(u_c)$ is not available. However, it is straightforward to demonstrate that  $v^*(u_c)$ is a
 continuous, monotone decreasing function of $u_c\in(0,1)$, with 
$v^*(u_c)\to 2^-$ as $u_c\to 0^{+}$
and $v^*(u_c)\to 0^+$ as $u_c\to 1^{-}$   \cite{Tisbury_etal2020a} , where $2$ is the minimum propagation speed of the PTW solution in the absence of 
cut-off ($u_c=0$).
Its asymptotic form can be approximated 
in the limits of $u_c\to 0^+$ and $u_c\to 1^-$. 
In the first limit, Brunet and Derrida \cite{BrunetDerrida1997} predicted
that  the value of $v^*(u_c)$  
is strongly influenced by the presence of a  cut-off, with 
\beq\elab{asympt1}
v^*(u_c)= 2-\frac{\pi^2}{(\ln u_c)^2}+O\left(\frac{1}{|\ln u_c|^3}\right)\quad\text{as $u_c\to 0^+$},
   \eeq
   so that in this limit, the difference between $v^*(u_c)$ and $2$ is only logarithmically small in $u_c$. 
This behaviour was rigorously verified by Dumortier, Popovic and Kaper
 \cite{Dumortier_etal2007} using geometric desingularization. 
 Higher order corrections were obtained by Tisbury and the authors  \cite{Tisbury_etal2020a}  using matched asymptotic expansions. The behaviour of $v^*(u_c)$  as $u_c\to 0^+$ is in contrast with the behaviour of  
$v^*(u_c)$ as $u_c\to 1^-$,  in which case it was established in  \cite{Tisbury_etal2020a}  that 
   \beq\elab{asympt2}
v^*(u_c)=(1-u_c)|f'(1)|^{1/2}+O((1-u_c)^2), 
\quad\text{as $u_c\to 1^{-}$},
   \eeq
   and so now vanishing algebraically in $(1-u_c)$.
There are no equivalent results for the spatially heterogenous version of (IBVP) ($A>0)$.

In this paper our aim is to understand how advection by a shear flow influences the shape of the  PTW solution and its propagation speed.  This understanding is important due {to} the fact that in 
a wide variety of environmental and engineering applications, associated chemical or biological reactive species
are transported by  fluid flows,  
with the simplest non-trivial flows being shear flows. 
In {this} case, it is natural to anticipate curved PTW solutions of the form 
$u(x,y,t)=U_T(x-vt,y)$ 
and propagation speed $v=\hat{v}(A,B,u_c)$ , where $\hat{v}(0,B,u_c)=v^*(u_c)$. The 
  functions $s:[0,1]\times\mathbb{R}^+\to \mathbb{R}$ and $U_T:\mathbb{R}\times[0,1]\to \mathbb{R}$ now satisfy
\beq
s(y,t)t^{-1}\to \hat{v}(A,B,u_c)
\quad\text{and}\quad
u(z+s(y,t),y,t)\to U_T(z,y),
\eeq
as $t\to\infty$ uniformly for $(z,y)\in\mathbb{R}\times[0,1]$ with $( s_t(y,t)-\hat{v}(A,B,u_c))=o(1)$
as $t\to\infty$ uniformly in $y\in[0,1]$. 
The question of 
existence of curved PTW solutions  for KPP reactions without cut-off (that is, $u_c=0$) and shear flows with bounded cross-sections and Neumann  boundary   conditions was 
considered in detail by Berestycki and Nirenberg \cite{BerestyckiNirenberg1992}.
They  
used the approach of sub- and super-solutions
on the associated semilinear elliptic boundary value problem  
to establish that 
a unique PTW solution exists 
for each propagation $v\geq v_m(A,B)$, 
where $v_m(A,B)$ is the (positive) minimum  available
propagation speed. 
This PTW solution is monotone decreasing in $z$ and has 
$\lim_{z\to\infty} U_T(z,y)=0$
and $\lim_{z\to-\infty} U_T(z,y)=1$, 
uniformly for $y\in[0,1]$. A number of works
\cite{Freidlin1985,MallordyRoquejoffre1995,Roquejoffre1997} have shown that in fact, for KPP reactions without cut-off,  the solution to the (IBVP)  approaches that  PTW  with minimum speed $v_m(A,B)$.

The elliptic boundary value problem whose solution provides the curved PTW    presents little immediate and direct information about the quantitative dependence of the propagation speed on
the parameters and the details of the profiles of the advective flow and reaction functions. 
A first characterization of the propagation speed was derived  by 
G\"artner and Freidlin \cite{GertnerFreidlin1979} for KPP reactions using probabilistic arguments.
This variational characterization expressed $v_m(A,B)$ in terms of  the principal eigenvalue of a certain linear eigenvalue problem that depends on $A$, $B$ and 
$\alpha(y)$. 
 Thus, $v_m(A,B)$ depends on the structure of $\alpha(y)$ but not the detailed structure  of $f(u)$ (other than it satisfying the KPP conditions). 
 An alternative characterization of  $v_m(A,B)$   was subsequently derived by
 Berestycki and Nirenberg \cite{BerestyckiNirenberg1992} in terms of a quadratic linear eigenvalue problem,
 where   $v_m(A,B)$ is determined by the requirement that the two principal eigenvalues are equal (the equivalence of the two characterizations is shown in \cite{Xin2000a}). 
 A variational characterization was developed in \cite{Heinze_etal2001} for the case of bistable and
 combustion type nonlinearity.

The eigenvalue problem determining  the  speed of propagation of the curved  PTW solution for KPP reaction functions without cut-off is readily solved numerically. An exact analytical solution is, however, not available  
even for simple shear flows. 
It can nevertheless be analysed in  asymptotic limits of the two parameters $A$ and $B$. 
For a channel width that is comparable to the advectionless front scale thickness,
i.e. when $B=O(1)$, it can be shown that a shear flow  enhances the speed of propagation, with the enhancement being monotonic in $A$ \cite{Berestycki2003}.
Further asymptotic results obtained in \cite{PapanicolaouXin1991} and 
\cite{Berestycki2003}  
respectively provide explicit expressions for  $v_m(A,B)$ with  
\beq\elab{vmB1}
v_m(A,B)=2+O(A^2/B) 
\quad\text{as}\quad A\to 0\quad\text{and}\quad
v_m(A,B)\sim A\alpha_M
\quad\text{as}\quad A\to \infty, 
\eeq 
  where the constant  $\alpha_M\equiv\sup_{y\in [0,1]}\alpha(y)>0$  
   was determined in  \cite{HaynesVanneste2014a}. 
For a channel width that is much smaller than the advectionless front scale thickness,
i.e. when $B\gg 1$,
$v_m(A,B)$ is, at leading order,  proportional to the square root of the   effective diffusivity  $\kappa_\text{eff}$
of the  advection--diffusion problem  \cite{Heinze_etal2001}
so that
\begin{subequations}\elab{vmBlarge}
\beq
v_m(A,B)\sim 2\sqrt{\kappa_\text{eff}}
\quad
\text{as}
\quad
B\to\infty, 
\eeq
uniformly for $A>0$. An explicit form for $\kappa_\text{eff}$
can be obtained using homogenization (see, for example, \cite{MajdaKramer1999, Camassa_etal2010}) which yields
\beq\elab{keff}
\kappa_\text{eff}=1+\left(A^2/B\right)\langle
\left(
\int_0^y\alpha(y')dy'
\right)^2
\rangle
\eeq
\end{subequations}
which holds uniformly for $A>0$. 
Here,  $\left\langle F (\cdot) \right\rangle\equiv \int_0^1  F(y) dy$. 
In the opposite limit, when the channel width  is much larger than the advectionless front scale thickness,
i.e. when $B\ll 1$,  
the PTW solutions are sharp-fronted and can be approximated by a single curve
where all the reaction takes place.  
Majda and Souganidis \cite{MajdaSouganidis1993} (see also \cite[Ch.\ 6]{Freidlin1985}) showed that 
a distinguished regime arises for  $A=O(1)$ 
in which case  the propagation speed satisfies  
\cite{Embid_etal1995,XinYu2013,TzellaVanneste2019} 
\beq\elab{vmBsmall}
v_m(A,B)\sim 2+A\alpha_M \quad\text{as}\quad B\to 0.
\eeq
A complete description of  $v_m(A,B)$ as $A/\sqrt{B}\to\infty$ can be readily deduced from the analysis in \cite{HaynesVanneste2014a}. 

A natural question is to what extent the above results apply in the presence of a cut-off in the reaction.

\subsection{Main results and paper structure}
In this paper we  consider (IBVP) in a number of natural asymptotic limits on the parameters $A$ and $B$, with $u_c\in (0,1)$ fixed.
One feature which we concentrate on, amongst others, is the existence of PTW solutions, and the detailed form of the associated propagation speed, 
$\hat{v}(A,B,u_c)$, which, for $u_c>0$, is unique. 
We use the theory of matched asymptotic expansions to establish the limiting form of (IBVP) and/or  
the associated PTW theory, and subsequently obtain expressions for $\hat{v}(A,B,u_c)$  
 and the boundary of the interface $\zeta(y)$ where $U_T(\zeta(y),y)=u_c$, in the three natural cases of $A\to \infty$, $A\to 0$ and $A=O(1)$. 
Our main conclusions are  that for $A\to\infty$ with $B=O(A^2)$, 
the solution to (IBVP) is, at leading order, described by an effective  equation devoid of advection with molecular diffusivity replaced by the effective diffusivity \eref{keff}. Thus,  
$\hat{v}(A,B,u_c)$ is at leading order enhanced by 
the shear flow by a factor proportional to the square root of $\kappa_{\text{eff}}$. 
For $A\to 0$  we consider {each of the complementary orderings $B^{-1}=o(1), B=O(1)$ and $B=o(1)$,} and find that $\hat{v}(A,B,u_c)$ is at leading order 
given by $v^*(u_c)$, the speed of propagation in the absence of a flow with 
 $\zeta(y)$ determined solely by the structure of the flow and the value of $v^*(u_c)$. 
The higher order correction to the speed is of $O(A)$ when $0<B \le O(A)$, initially increasing from zero as $B^{1/2}$, until it achieves a maximum  point located at a value $B=B_M(A) = O(A)$. Thereafter it is decreasing with increasing $B$, becoming of $O(A^2)$ for $B\ge O(1)$, and decreasing at a rate of $B^{-1}$.
At  the same time, the  structure of the interface   becomes increasingly {deformed}.  
Finally, for $A=O(1)$ with $B\to\infty$, the solution to (IBVP) is, at leading order, effectively described by the (IBVP) obtained for $A=B=0$. 
The situation is different for $A=O(1)$ with $B\to 0$, when the leading order term of the propagation speed   also includes $A\alpha_M$, 
where $\alpha_M$ corresponds to the maximum velocity in the channel, whilst the correction is of $O(B^{\frac{1}{2}})$. 
In this case, the interface is most {deformed}. 
Contrasting our results against \eref{vmB1}, \eref{vmBlarge} and \eref{vmBsmall}, \emph{we conclude that the effect of a shear flow on the speed of propagation of the PTW solution is similar  with and without cut-off in the KPP reaction, as may be anticipated}. Finally, it should be noted that in each of the overlapping limits within and between sections 3--5, it is readily verified that the respective asymptotic expressions for $\hat{v}(A,B,u_c)$ match together according to the classical Van Dyke asymptotic matching principle \cite{VanDyke1975} {(see also \cite[Ch.\ 8]{Miller} for an applied analysis point of view)}.

The paper is structured as follows. Section 2 focuses on reformulating (IBVP) as an equivalent moving boundary evolution problem that we refer to as (QIVP), and as a preliminary, we examine the structure of the solution to (QIVP) at small times. We then move on to describe two equivalent elliptic boundary value problems that govern the existence of possible PTW solutions to (IBVP) and (QIVP) respectively, and explore their structure as $u_c\to 1$. 
Sections 3, 4 and 5 are devoted to each of the three cases of $A\to \infty$, $A\to0$ and $A=O(1)$ respectively.
Throughout, we  illustrate the theory by evaluating expressions (either directly or through numerical integration) for  two classical shear flows: the plane Couette flow given by the linear shear
\begin{subequations}\elab{CouettePoiseuille}
\begin{align}
\alpha(y)&=y-1/2,  \elab{Couette} \\
\intertext{and the plane Poiseuille flow given by}
\alpha(y)&=-2y^2+2y-1/3,
\elab{Poiseuille}
\end{align}
\end{subequations}
with the constants fixed by the requirement that both of these shear flows have mean zero flow, their half-channel mean strain rate is the same and the maximum of the Poiseuille flow is located at the channel centre line $y=1/2$.
The paper ends with the concluding section 6.

\section{Problem formulation}
We begin this section by reformulating the parabolic evolution problem (IBVP) as an equivalent moving boundary evolution problem that we refer to as (QIVP). This reformulation will prove convenient at a number of stages throughout the paper. We then move on to describe two equivalent elliptic boundary value problems derived from (IBVP) and (QIVP) whose solution provides  permanent form travelling waves. We then consider some general results relating to this evolution problem, and its reduction for permanent form travelling waves. 

\subsection{The moving boundary evolution problem}
In this subsection we begin by developing a {modification and then a reformulation} of (IBVP). Due to the discontinuity in $f_c(u)$ when $u=u_c$, there will be a lack of full regularity in classical solutions to (IBVP) {(this is readily established, employing an argument by contradiction, after it is shown, via the maximum principle, that any fully classical solution to (IBVP) must have $u_x<0$ on $\mathbb{R}\times (0,T]$, for any $T>0$, and so, at each $t\in (0,T]$, there exists a unique smooth curve $x=x(y,t)$ upon which $u=u_c$. The contradiction then follows by choosing an interior point on this curve, and then examining the limits in 
 \eref{KPP} as this point is approached from left and right, under fully classical regularity and with reaction function \eref{BDreaction})}. 
{Therefore, for a solution to exist at all to (IBVP), we must mildly relax regularity and the notion of fully classical solution.  
Specifically, across local, spatial level curves along which $u=u_c$, we 
admit the possibility of a} \emph{jump discontinuity in the second spatial partial derivative of $u$ in the direction normal to such a curve}.
{Furthermore, from \cite{Tisbury_etal2020a}, we anticipate that such a spatial curve will be unique. 
With this in mind, we adopt
a reformulation of (IBVP), \emph{which under the above  restrictions}, can be written as an equivalent moving boundary evolution problem. 
} 
First, as a measure of convenience in this reformulation, we introduce a simple rescaling of coordinates via,  
\begin{equation*} 
(x,y)=(x',L^{-1}y')
\end{equation*}
and the problem domain now has
$(x',y',t)\in \mathbb{R}\times[0,L]\times \overline{\mathbb{R}}^+=\bar{D}$
with $D=\mathbb{R}\times(0,L)\times \mathbb{R}^+$ and $L=1/\sqrt{B}.$
We next introduce the smooth interface in $D$ as
 \beq\elab{eqn2.3}
\mathcal{L}_{IBVP} = \{(x',y',t)\in D : x'=s(t) + \zeta(y',t)\}
 \eeq
where 
\begin{equation}\elab{regularity_int}
	s\in C^1(\mathbb{R}^+)\cap C(\bar{\mathbb{R}}^+)\quad\text{and}\quad\zeta\in C^{1,1}([0,L]\times{\mathbb{R}}^+) \cap C([0,L]\times\bar{\mathbb{R}}^+), 
\end{equation}	
	and are chosen so that
\beq\elab{eqn2.4}
u(s(t)+\zeta(y',t),y',t) = u_c~~\forall~~(y',t)\in[0,L]\times{\mathbb{R}}^+,
\quad\text{with}\quad
\int_0^L{\zeta(y',t)}dy'=0~~\forall~~t\in\mathbb{R}^+.
\eeq
At each $t>0$, we observe that $\mathcal{L}_{IBVP}$ represents the spatial curve for $(x',y')\in \mathbb{R}\times[0,L]$ on which $u=u_c$. In association with $\mathcal{L}_{IBVP}$, we introduce the regions
\begin{subequations}
\beq\elab{eqn2.6}
D_R = \{(x',y',t) : (y',t)\in(0,L)\times{\mathbb{R}^+}, x'\in(s(t) + \zeta(y',t),\infty)\}
\eeq
and
\beq\elab{eqn2.7}
D_L = \{(x',y',t) : (y',t)\in(0,L)\times\mathbb{R}^+, x'\in(-\infty,s(t)+\zeta(y',t)\}
\eeq
\end{subequations}
with $u\ge u_c$ in $\bar{D}_L$ and $u\le u_c$ in $\bar{D}_R$. In this context, a classical solution will have $u:\mathbb{R}\times[0,L]\times\overline{\mathbb{R}}^+\to\mathbb{R}$ having regularity 
\beq\elab{eqn2.8}
u \in BC((\mathbb{R}\times[0,L]\times\overline{\mathbb{R}}^+)\setminus(\{0\}\times[0,L]\times\{0\}))~ \cap~ C^{1,1,1}(\mathbb{R}\times[0,L]\times \mathbb{R}^+)~\cap~ C^{2,2,1}(D_L\cup D_R).
\eeq
{Here $BC$ represents bounded and continuous whilst, for example, $C^{1,1,1}$ represents continuous with continuous partial derivatives in the first three arguments.} The moving boundary problem is then readily formulated as follows. First, it is convenient to introduce the mean shifted coordinate along the channel, namely,
\[ 
\xi'=x'-s(t)~~\forall~~(x',t)\in\mathbb{R}\times\overline{\mathbb{R}}^+.
\] 
  with the original interface now located at 
 \beq\elab{eqn2.23}
\mathcal{L}_{QIVP} = \{(\xi',y',t):\xi'=\zeta(y',t),~(y',t)\in (0,L)\times \mathbb{R}^+\}.
\eeq
At each $t> 0$, we note that a spatial normal to $\mathcal{L}_{QIVP}$, pointing left to right, is given by
\begin{equation*} 
\textbf{n}(\zeta(y',t),y')=(1,-\zeta_{y'}(y',t))~~\forall~~(y',t)\in (0,L)\times \mathbb{R}^+.
\end{equation*} 
\begin{taggedsubequations}{QIVP}\label{QIVP}
\beq\elab{eqn2.10}
u_t + (A\bar{\alpha}(y') - \dot{s}(t))u_{\xi'} = \nabla^2u + f_c(u),~~(\xi',y',t)\in Q_L\cup Q_R,
\eeq
\beq\elab{eqn2.11}
u\ge u_c~\text{in}~Q_L,~~u\le u_c~\text{in}~Q_R,
\eeq
\beq\elab{eqn2.12}
u(\xi',y',0)=H(-\xi'),~~(\xi',y')\in \mathbb{R}\times[0,L],
\eeq
\beq\elab{eqn2.13}
u(\xi',y',t)\to
\begin{cases}
1~~\text{as}~~\xi'\to -\infty,\\
0~~\text{as}~~\xi'\to \infty,
\end{cases}
\eeq
uniformly for $(y',t)\in [0,L]\times[0,T]$, any $T>0$,
\beq\elab{eqn2.14}
u_{y'}(\xi',0,t)=u_{y'}(\xi',L,t)=0,~~(\xi',t)\in \mathbb{R}\times\overline{\mathbb{R}}^+,
\eeq
and at the moving boundary, {continuity of $u$ and its two spatial partial derivatives requires (with a direct calculation of directional derivatives),}
\beq\elab{eqn2.15}
u(\zeta(y',t)^+,y',t)=u(\zeta(y',t)^-,y',t)=u_c,~~(y',t)\in [0,L]\times\mathbb{R}^+,
\eeq
\beq\elab{eqn2.16}
s(0)=0,~\zeta(y',0)=0,~~y'\in [0,L],
\eeq
\beq\elab{eqn2.17}
[\nabla u\cdot \bm{n}(\zeta(y',t),y',t))]_L^R=0,~~(y',t)\in (0,L)\times\mathbb{R}^+.
\eeq
In the above, $\nabla=(\partial_{\xi'},\partial_{y'})$ is the usual gradient operator in the coordinates $(\xi',y')$, and,
\beq\elab{eqn2.18}
\begin{split}
Q_R & = \{(\xi',y',t) : (y',t)\in (0,L)\times \mathbb{R}^+,~\xi'\in (\zeta(y',t),\infty)\},\\
Q_L &= \{(\xi',y',t) : (y',t)\in (0,L)\times \mathbb{R}^+,~\xi'\in (-\infty,\zeta(y',t))\},
\end{split}
\eeq
with $[\cdot]_L^R$ indicating the difference in limits across $\mathcal{L}_{QIVP}$ from $Q_L$ to $Q_R$. In addition, $\bar{\alpha}(y')=\alpha(y'/L)~~\forall~~y'\in [0,L]$
with, from  \eref{eqn1.4} and \eref{eqn2.4},
\beq\elab{eqn2.21}
\int_0^L{\bar{\alpha}(y')}dy' = 0,
\quad\text{and}\quad
\int_0^L{\zeta(y',t)}dy'=0~~\forall~~t\in \mathbb{R}^+.
\eeq
\end{taggedsubequations}
The regularity condition given in \eref{eqn2.8} becomes,
\beq\elab{eqn2.24}
u \in BC((\mathbb{R}\times[0,L]\times[0,\infty))\setminus(\{0\}\times[0,L]\times\{0\}))~ \cap~ C^{1,1,1}(\mathbb{R}\times[0,L]\times(0,\infty))~\cap~ C^{2,2,1}(Q_L\cup Q_R).
\eeq
It is useful to observe that on using the regularity requirements on $u$ as given in \eref{eqn2.24}, together with the earlier regularity requirements on $\zeta$ and $s$ \eref{regularity_int},  
we can readily establish that $u_{\xi'\xi'}$, $u_{y'y'}$ and $u_{\xi'y'}$ each have finite limits from $Q_L$ and $Q_R$ as each point on $\mathcal{L}_{QIVP}$ is approached. In addition, on using the translation invariance of (QIVP) in the coordinate $\xi'$, together with application of the classical parabolic strong maximum principle and comparison theorem (in both $Q_L$ and $Q_R$ respectively) we are able to establish the following basic properties for (QIVP), namely,
\beq\elab{eqn2.26}
0<u(\xi',y',t)<u_c~~\text{in}~~Q_R,
\quad
u_c<u(\xi',y',t)<1~~\text{in}~~Q_L,~\forall t>0,
\eeq
and,
\beq\elab{eqn2.28}
u(\xi',y',t)~~ \text{is strictly monotone decreasing with}~~\xi'\in \mathbb{R}~~ \text{for fixed}~(y',t)\in (0,L)\times(0,\infty).
\eeq
Use of the regularity conditions on $u$ and $\zeta$, together with the moving boundary conditions \eref{eqn2.15} and \eref{eqn2.17}, the chain rule, and the PDE 
\eref{eqn2.10}, in both $Q_L$ and $Q_R$, it is also readily established
{(using continuity of the two spatial partial derivatives  $u_{ss}$ and $u_{ns}$ across the interface, together with 
\eref{eqn2.10}, when written in terms of tangential-normal coordinates $(s,n)$ at the interface $\mathcal{L}_{QIVP}$)} that
\beq\elab{eqn2.29}
[ u_{nn}]_L^R = f_c, 
\eeq
at all points on $\mathcal{L}_{QIVP}$, with $u_{nn}$ being the second spatial derivative of $u$ with respect to distance in the direction of $\textbf{n},$ as the point on $\mathcal{L}_{QIVP}$ is approached.

 {We are now able to construct, via the method of matched asymptotic expansions, a solution to (QIVP) as $t\to 0^+$ which gives a uniform approximation to $u$ throughout the spatial domain.
 In particular, this analysis provides information regarding the early dynamics of the interface in (QIVP), and the relative importance of advection, diffusion and reaction in the dynamics. 
 We note that the ability to construct such an approximation gives a strong indication that (QIVP) is, at least locally in $t$, well-posed.  
The details of this construction are extensive and as such are relegated to  the Appendix. }

Finally we consider, at a given $t=O(1)^+$, a local analysis close to the point of intersection of the spatial  moving boundary and the channel boundary, at $(\xi',y')=(\zeta(0,t),0)$\footnote{{In this context, for any variable $\lambda$, we will henceforth write $\lambda=O(1)>0$ as $\lambda=O(1)^+$, and correspondingly, $\lambda=O(1)<0$ as $\lambda=O(1)^-$}.}.  
 After some calculation, it follows from \eref{eqn2.10}, \eref{eqn2.11},\eref{eqn2.14} and \eref{eqn2.15} that,
\beq\elab{eqn2.30}
u(r,\theta,t) = 
\begin{cases}
u_c - A_+(t)r^{\frac{\pi}{2\alpha(t)}}\text{cos}(\frac{\pi\theta}{2\alpha(t)}) + o(r^{\frac{\pi}{2\alpha(t)}})~~\text{as}~~r\to 0^+ ~~\text{in}~~ \bar{Q}_+,\\
u_c + A_-(t)r^{\frac{\pi}{2(\pi-\alpha(t))}}\text{cos}(\frac{\pi(\pi-\theta)}{2(\pi-\alpha(t))}) + o(r^{\frac{\pi}{2(\pi-\alpha(t))}})~~\text{as}~~r\to 0^+ ~~\text{in}~~ \bar{Q}_-,
\end{cases}
\eeq
Here $A_+(t)$ and $A_-(t)$ are positive (via \eref{eqn2.26}), smooth, globally determined functions for $t=O(1)^+$, whilst $(r,\theta)$ are plane polar coordinates, given by
\begin{subequations}
\beq\elab{eqn2.31}
\xi'=\zeta(0,t)+r\text{cos}\theta,~~y'=r\text{sin}\theta,
\eeq
with $r\ge0$ and $0\le \theta \le \pi$. In addition,
\beq\elab{eqn2.32}
\text{tan}\alpha(t) = \zeta_{y'}(0,t)~~t>0,
\eeq
\end{subequations}
so that, locally,
\beq\elab{eqn2.33}
\bar{Q}_+ = \{(r,\theta):r\ge 0,~~0\le \theta\le \alpha(t)\},
\quad
\bar{Q}_- = \{(r,\theta):r\ge 0,~~\alpha(t)\le \theta\le \pi\},
\eeq
with
\beq\elab{eqn2.35}
0<\alpha(t)<\pi
\eeq
at each $t=O(1)^+$. It remains to apply condition \eref{eqn2.17}, which finally requires that,
\beq\elab{eqn2.36}
\alpha(t)=\frac{1}{2}\pi~~\text{and}~~A_+(t)=A_-(t).
\eeq
 A similar analysis follows close to the intersection point of the spatial moving boundary and the channel boundary at $(\xi',y')=(\zeta(L,t),L).$ We observe from \eref{eqn2.36} that the spatial moving boundary makes normal contact with the channel boundary at each $t>0$. It should be noted that the initial condition \eref{eqn2.12} requires that $\zeta(0,t)\to 0$ as $t\to 0^+$. 
 In addition, the structure of the solution to (QIVP) when $t$ is small, which is 
 developed in the {Appendix, does in fact allow us to obtain some additional information regarding the indeterminate functions $A_{\pm}(t)$. In particular, asymptotic matching with particular boundary regions, which are located, respectively when $y'=O(t^{1/2})$ and $y'=L-O(t^{1/2})$ 
  (labelled as  regions $\mathrm{I}_{L,R}$ in the Appendix)  determines (via Van Dyke's matching principle)}   that
 \beq\elab{eqn2.37'}
 A_{\pm}(t) = O(t^{-\frac{1}{2}})~~\text{as}~~t\to 0^+.
 \eeq
 This condition, together with \eref{eqn2.36}, determines that \emph{when $t$ is small}, the spatial region of validity of the expansion \eref{eqn2.30}, local to the free boundary contact point with the wall at $y'=0$, requires the restriction $0 \le r \ll t^{\frac{1}{2}}$ with $0 \le \theta \le \pi$.

We next consider permanent form travelling front structures which may develop as large-$t$ attractors in the solution to (IBVP), and equivalently (QIVP).

\subsection{Permanent form travelling fronts}
We anticipate that as $t\to \infty$ a permanent form travelling front structure will develop in the solution to (IBVP) (and equivalently (QIVP)), advancing with a non-negative, constant propagation speed, and allowing for the transition from the unreacted state $u=0$ to the fully reacted state $u=1$. {This anticipation is supported by the theory developed in the earlier papers by Tisbury and the authors, \cite{Tisbury_etal2020a} and \cite{Tisbury_etal2020b}, in the case when the shear flow is absent. With this in mind,} we begin by formulating the elliptic boundary value problems (which are equivalent) associated with permanent form travelling front solutions in both the formulations (IBVP) and (QIVP). {As will be seen in the later sections, one or other of these equivalent formulations usually presents itself as the most natural to work with in each limiting case that we consider.} In what follows, we will refer to a permanent form travelling front solution to either of (IBVP) or (QIVP) as a PTW solution.

A PTW solution, with propagation speed $v\ge0$, is a non-negative solution to (IBVP),  
which depends only upon $y$ and the travelling coordinate,
\[
z=x-vt.
\]
Thus, $U_T:\mathbb{R}\times[0,1]\to \mathbb{R}$ is a PTW soluton to (IBVP), with propagation speed $v\ge0$, when,
\begin{taggedsubequations}{BVP}\label{BVP}
\beq\elab{eqn2.53}
U_{Tzz} + BU_{Tyy} + (v-A\alpha(y))U_{Tz} + f_c(U_T) = 0,~~(z,y)\in \mathbb{R}\times(0,1), 
\eeq
with $U_T\in C^1(\mathbb{R}\times[0,1])\cap C^2(\mathbb{R}\times((0,1)\setminus \mathcal{L}_{BVP}))$, {where the interface is now,}
\beq\elab{eqn2.54}
\mathcal{L}_{BVP}= \{(z,y)\in \mathbb{R}\times(0,1): u(z,y)=u_c\}.
\eeq
In addition,
\beq\elab{eqn2.55}
U_T(z,y)\ge0~~\forall~~(z,y)\in \mathbb{R}\times[0,1],
\eeq
\beq\elab{eqn2.56}
U_{Ty}(z,0)=U_{Ty}(z,1)=0~~\forall~~z\in \mathbb{R},
\eeq
\beq\elab{eqn2.57}
U_T(z,y)\to
\begin{cases}
0~~\text{as}~~z\to \infty, \\
1~~\text{as}~~z\to {-\infty}
\end{cases}
\eeq
\end{taggedsubequations}
uniformly for $y\in [0,1]$. This elliptic boundary value problem, {which we henceforth refer to as (BVP)}, may be regarded as a nonlinear eigenvalue problem, with eigenvalue $v\ge0$. 
Now let $U_T:\mathbb{R}\times[0,1]\to \mathbb{R}$ be a PTW with propagation speed $v\ge 0$. We observe that, on any closed bounded interval $I\subset \mathbb{R}$, then,
\beq\elab{2.58}
f_c(\lambda)\lambda^{-1}~\text{is bounded}~~\forall~~\lambda\in I,
\eeq
via \eref{BDreaction} and \eref{KPPreaction}. It then follows from \eref{eqn2.53}, \eref{eqn2.56} and \eref{eqn2.57}, together with the strong elliptic maximum principle (see, for example,  Ch.\ 3 in 
\cite{GilbargTrudinger2001}), that,
\beq\elab{eqn2.59}
U_T(z,y) > 0~~\forall~~(z,y)\in \mathbb{R}\times[0,1].
\eeq
A similar argument, using \eref{eqn2.53}, \eref{eqn2.56} and \eref{eqn2.57}, with \eref{BDreaction}, \eref{KPPreaction}, and the strong elliptic maximum principle establishes that
\beq\elab{eqn2.60}
U_T(z,y) < 1~~\forall~~(z,y)\in \mathbb{R}\times[0,1].
\eeq
With the moving boundary reformulation in mind, we {return to the  formulation of (QIVP) which determines} that $\mathcal{L}_{BVP}$ can be represented in the following form,
\beq\elab{eqn2.61}
\mathcal{L}_{BVP}= \{(z,y) : z=\gamma(y)~~ \text{for}~~ y\in [0,1]~~ \text{with}~~ U_T(\gamma(y),y)=u_c\}
\eeq 
with (via the regularity required on $U_T:\mathbb{R}\times[0,1]\to \mathbb{R}$) $\gamma\in C^1([0,1])$, and $u\ge u_c$ to the left of $\mathcal{L}_{BVP}$ and $u\le u_c$ to the right of $\mathcal{L}_{BVP}$, respectively. With this structure, we can then apply the strong maximum and minimum principle, in each of the regions to the right and left of $\mathcal{L}_{BVP}$, to establish that,
\beq\elab{eqn2.62}
U_{Tz}(z,y)<0~~\forall~~(z,y)\in \mathbb{R}\times(0,1)
\quad\text{with}\quad
U_T(z,y)
\begin{cases}
>u_c~~\text{to left of}~~\mathcal{L}_{BVP},\\
<u_c~~\text{to right of}~~\mathcal{L}_{BVP}.
\end{cases}
\eeq
The problem (BVP) is translation invariant in the coordinate $z$, and we will fix this invariance by henceforth requiring that
\beq\elab{eqn2.64}
\int_0^1{\gamma(y)}dy = 0.
\eeq
We now consider the equivalent moving boundary formulation for a PTW solution. With translational invariance fixed through \eref{eqn2.64}, then a PTW is simply a \emph{steady solution} of (QIVP), and so the  spatial moving boundary in (QIVP) for a PTW, with propagation speed $v\ge 0$, is now located at,
\beq\elab{eqn2.65}
\mathcal{L}_{QBVP}= \{(\xi',y'):~~\xi'=\zeta(y'),~~y'\in [0,L]\}
\eeq
where, as the solution is now steady, the dependence on $t$ has been omitted, whilst in relation to (BVP),
\beq\elab{eqn2.66}
\xi'=z,~~\zeta(y')=\gamma(y'L^{-1})
\eeq
for each $y'\in [0,L]$, with $y$ and $y'$ related as before, whilst translational invariance \eref{eqn2.64} becomes (following \eref{eqn2.4} and \eref{eqn2.21}) 
\beq\elab{eqn2.67}
\int_0^L{\zeta(y')}dy' = 0
\eeq
with
\beq\elab{eqn2.68}
\dot{s}(t) = v.
\eeq
{The (BVP) formulation becomes the equivalent moving boundary formulation that we henceforth refer to as (QBVP):} 
\begin{taggedsubequations}{QBVP}\label{QBVP}
\beq\elab{eqn2.69}
  \nabla^2U_T - (A\bar{\alpha}(y') - v)U_{T\xi'} + f_c(U_T) = 0,~~(\xi',y')\in Q_L\cup Q_R,
\eeq
with,
\beq\elab{eqn2.70}
U_T\ge u_c~\text{in}~Q_L,~~U_T\le u_c~\text{in}~Q_R,
\eeq
and
\beq\elab{eqn2.71}
U_T(\xi',y')\to
\begin{cases}
1~~\text{as}~~\xi'\to -\infty,\\
0~~\text{as}~~\xi'\to \infty,
\end{cases}
\eeq
uniformly for $y'\in [0,L]$, whilst
\beq\elab{eqn2.72}
U_{Ty'}(\xi',0)=U_{Ty'}(\xi',L)=0,~~\xi'\in \mathbb{R}.
\eeq
At the boundary $\mathcal{L}_{QBVP}$,
\beq\elab{eqn2.73}
U_T(\zeta(y')^+,y')=U_T(\zeta(y')^-,y')=u_c,~~y'\in [0,L],
\eeq
and
\beq\elab{eqn2.74}
[\nabla U_T\cdot \bm{n}(\zeta(y'),y'))]_L^R=0,~~y'\in (0,L).
\eeq 
\end{taggedsubequations}
Throughout the rest of the paper, we will consider PTW solutions via either (BVP) or (QBVP), depending upon convenience. Before moving on to detailed analysis, we first make some preliminary observations regarding PTW solutions. Following the theory developed in \cite{Tisbury_etal2020a}, for the case when $\bar{\alpha}(y')=0$ for $y'\in [0,L]$, we anticipate that, for a \emph{given} shear flow profile $\bar{\alpha}(y')$, then (QBVP) (and equivalently (BVP)) has a PTW solution if and only if
\beq\elab{eqn2.75}
v = \hat{v}(A,B,u_c)>0 
\eeq
with $\hat{v}:(0,\infty)\times \mathbb{R}^+\times(0,1)\to \mathbb{R}$ 
a continuous function, and that this PTW solution is unique. We will see that this is supported by  the asymptotic results developed in the following sections. In addition, following \cite{Tisbury_etal2020a}, we anticipate that,
\beq\elab{eqn2.76}
\hat{v}(A,B,u_c) \to 0~~\text{as}~~u_c\to 1
\eeq
whilst,
\beq\elab{eqn2.77}
\hat{v}(A,B,u_c) \to v_m(A,B)~~\text{as}~~u_c\to 0
\eeq
where, following \cite{BerestyckiNirenberg1992} $v_m(A,B)>0$ is the \emph{minimum propagation speed for PTW solutions to the corresponding KPP problem in the absence of cut-off}. 

\subsection{Asymptotic structure of the solution to (QBVP) as $\bm{u_c\to 1}$.}
To complete this section we consider in more detail the form of $\hat{v}(A,B,u_c)$ as $u_c\to 1$ with $A,B=O(1)$. For this purpose, it is convenient to use the moving boundary formulation for PTW solutions, namely (QBVP).
{We integrate  equation \eref{eqn2.69} over $Q_L$ and $Q_R$ and apply Green's Theorem in the plane,   using conditions \eref{eqn2.21}, \eref{eqn2.70} -- \eref{eqn2.74}}
 and taking care with the existence of the associated improper integrals, which gives, in general,
\beq\elab{eqn2.78}
\hat{v}(A,B,u_c) = \frac{1}{L} \iint_{\bar{Q}_L}\ f_c(U_T(\xi',y')) \,d\xi' \,dy'.
\eeq
{We may now restrict attention to the spatial subdomain $Q_L$. In $Q_L$ it follows directly from \eref{eqn2.15} and \eref{eqn2.28} that $U_T = 1- O((1-u_c))$ whilst a balancing of terms in
the PDE \eref{eqn2.10}
in $Q_L$ determines that $\hat{v}=O((1-u_c))$, as $u_c\to 1^-$. Therefore we expand in $Q_L$ as,}
\beq\elab{eqn2.79}
\hat{v}(A,B,u_c) = (1-u_c) \bar{v}(A,B) + o((1-u_c)),
\eeq
\beq\elab{eqn2.80}
U_T(\xi',y') = 1 - (1-u_c)\bar{u}(\xi',y') + o((1-u_c)),
\eeq
with $(\xi',y')\in \bar{Q}_L$, and $A,B = O(1)$. After some straightforward calculation, we find that,
\beq\elab{eqn2.81}
\bar{u}(\xi',y') = c\tilde{\phi}(y')e^{\lambda_0\xi'},~~(\xi',y')\in \bar{Q}_L,
\eeq
with $\lambda=\lambda_0>0$ being the smallest positive eigenvalue \emph{with a principal eigenfunction}, for the quadratic linear eigenvalue problem, 
\beq\elab{eqn2.82}
\tilde{\phi}'' + (\lambda^2 - A{\bar\alpha}(y')\lambda + f_c'(1))\tilde{\phi} = 0,~~y'\in (0,L),
\eeq
\beq\elab{eqn2.83}
\tilde{\phi}'(0) = \tilde{\phi}'(L) = 0,
\eeq
\beq\elab{eqn2.85'}
\tilde{\phi}(y')>0~~\forall~~y'\in [0,L]
\eeq
with normalisation, 
\beq\elab{eqn2.85}
\int_{0}^{L}{\tilde{\phi}(y')}dy' = 1,
\eeq
and where $c$ is a constant to be determined (it is readily established that, with $f_c'(1) <0,$ this quadratic linear eigenvalue problem has exactly two distinct eigenvalues which have principal eigenfunctions, and these two eigenvalues are real, and have opposite sign).
We now determine $\zeta(y')$ using \eref{eqn2.80} and \eref{eqn2.81} in the remaining boundary condition \eref{eqn2.73} to obtain,
\beq\elab{eqn2.86}
\zeta(y') = - \frac{1}{\lambda_0}\log{(c\tilde{\phi}(y'))}
\eeq
for $y'\in [0,L]$. It now remains to apply the final condition \eref{eqn2.67} after which we arrive at,
\beq\elab{eqn2.87}
c = e^{-M}
\quad
\text{with}
\quad
M = \langle\log{\tilde{\phi}}\rangle = \frac{1}{L}\int_{0}^{L}{\log{\tilde{\phi}(s)}}ds.
\eeq
Thus, the moving boundary is located at,
\beq\elab{eqn2.89}
\xi' = \frac{1}{\lambda_0}(\langle\log{\tilde{\phi}}\rangle - \log{\tilde{\phi}(y')}) + o(1),
\eeq
as $u_c\to 1^-$ with $y'\in [0,L]$. Also, from \eref{eqn2.78} -- \eref{eqn2.81}, with \eref{eqn2.87}, we obtain
\beq
\bar{v}(A,B) = \frac{1}{L}c|f_c'(1)| \iint_{\bar{Q}_L}\ \tilde{\phi}(y')e^{\lambda_0\xi'}  \,d\xi' \,dy',
\eeq
recalling that both $\lambda_0$ and $\tilde{\phi}$ depend upon $A$ and $B$. The double integral can be evaluated directly, using \eref{eqn2.87} and \eref{eqn2.89}, as having value $\frac{L}{\lambda_0c}$. Thus we finally obtain,
\beq\elab{smallucasymp}
\bar{v}(A,B) = \frac{|f_c'(1)|}{\lambda_0(A,B)}.
\eeq
Clearly, in the absence of advection (corresponding to $A=0$) a direct calculation gives $\lambda_0(0,B)=\sqrt{|f_c'(1)|}$. In this case, equation \eref{smallucasymp} recovers the asymptotic result obtained in \cite{Tisbury_etal2020a}. In addition, it is readily established that
\beq\elab{eqn2.90}
\lambda_0(A,B) = H(I(A,B))~~\forall~~(A,B)\in \overline{\mathbb{R}}^+ \times \mathbb{R}^+
\eeq
where $H:\mathbb{R}\to \mathbb{R}$ is the \emph{positive-valued and strictly monotone increasing} function
\[
H(X) = \frac{1}{2}\left(X + \sqrt{(X^2 +4|f'_c(1)|)}\right)~~\forall~~X\in \mathbb{R}
\]
and
\[
I(A,B) = A\int_0^L{\overline{\alpha}(s)\overline{\phi}(s)}ds~~\forall~~(A,B)\in \overline{\mathbb{R}}^+ \times \mathbb{R}^+
\]
recalling that {$B=L^{-2}$}. We now set
\[
m(B) = \text{inf}_{s\in[0,L]}\overline{\alpha}(s)<0
\quad
\text{and}
\quad
M(B) = \text{sup}_{s\in[0,L]}\overline{\alpha}(s)>0.
\]
It then follows from \eref{eqn2.85'} and \eref{eqn2.85} that
\[
Am(B) < I(A,B) < AM(B)~~\forall~~(A,B)\in \overline{\mathbb{R}}^+ \times \mathbb{R}^+
\]
after which we obtain, via \eref{eqn2.90} and monotonicity of $H(\cdot)$, that
\beq\elab{bds}
H(Am(B)) < \lambda_0(A,B) < H(AM(B))~~\forall~~(A,B)\in \overline{\mathbb{R}}^+ \times \mathbb{R}^+.
\eeq
We observe that,
\[
 H(Am(B)) < \lambda_0(0,B)=\sqrt{|f'_c(1)|} < H(AM(B))~~\forall~~(A,B)\in \mathbb{R}^+ \times \mathbb{R}^+.
\]
Thus, although the bounds \eref{bds} constrain (via \eref{smallucasymp}) the   propagation speed as $u_c\to 1^-$,  they do not indicate whether the propagation speed is enhanced or retarded by the inclusion of advection. 

\begin{figure}
		 \begin{center}
 	     \includegraphics[width=0.49\textwidth]{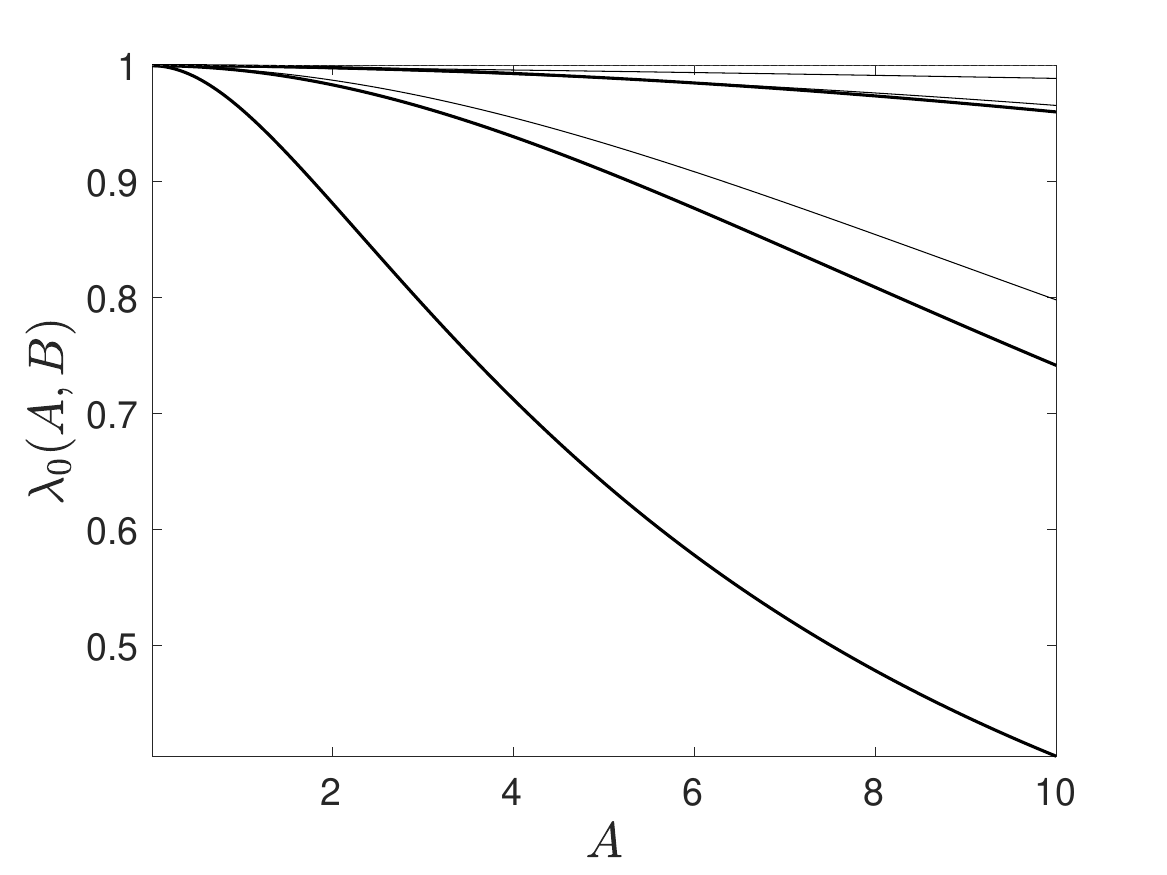}
 \includegraphics[width=0.49\textwidth]{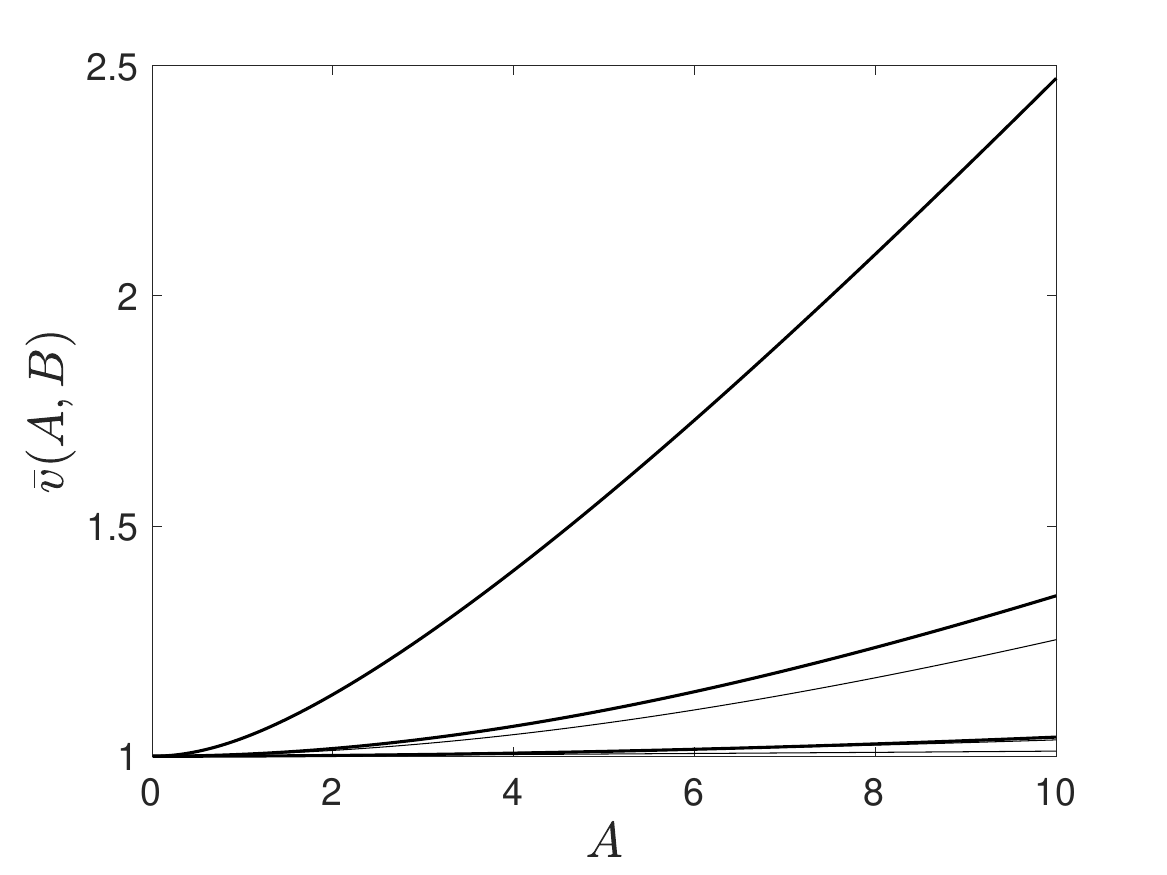}
		  \end{center}
    \caption {A graph of (left) the principal eigenvalue $\lambda_0(A, B)$ to 
    the quadratic linear eigenvalue problem
    \eref{eqn2.82}--\eref{eqn2.83} 
    and  (right) $\bar v(A,B)$  determined for $f'_c(1)=-1$ using \eref{smallucasymp}. 
     Results are obtained numerically for
     the Couette     (thick solid lines) and the 
     Poiseuille   flow (thin solid lines) for (left) $B=0.1$ (bottom), $1$ (middle)
     and $10$ (top) and {in reverse order} for (right). }
 
  \flab{largecutoff}  
  \end{figure}

We now  compute the positive principal eigenvalue 
$\lambda_0(A, B)$ to the quadratic linear eigenvalue problem for the Couette flow \eref{Couette} and the Poiseuille flow  \eref{Poiseuille} and use \eref{smallucasymp} with $f'_c(1)=-1$ to evaluate $\bar v(A,B)$ in these two example cases. We use 
a standard second-order finite-difference discretization of \eref{eqn2.82}--\eref{eqn2.83}. The resulting matrix eigenvalue problem is solved for a range of values of $A$ and $B$ using MATLAB’s routine \texttt{polyeigs}. We choose the spatial resolution $h$ to satisfy $1/h= 200$. The results are shown in Figure \fref{largecutoff}. It is clear that for these particular flows, the propagation speed is enhanced by advection, increasing with $A$ and decreasing with $B$. Furthermore, for all values of $A$ and $B$ considered, the propagation speed  in the Couette flow remains larger   than in the plane Poiseuille flow.

We   now proceed to consider (IBVP) or, equivalently (QIVP) and (BVP) or, equivalently (QBVP), in a number of 
significant asymptotic limits in the two key parameters $(A,B)$.

\section{Slowly varying front with strong advection: \textbf{$B^{-1}=o(1)$} as \textbf{$A\to \infty$}}
In this section we consider (IBVP) when the advection velocity scale is large compared to the advectionless front propagation speed, and the channel width is small relative to the advectionless front scale thickness. In terms of the two parameters $A$ and $B$, we can formalise this by considering $A\gg1$ and $B\gg1$, with the principal balances in the PDE \eref{KPP} leading us to consider the limit $A\to \infty$, with $B=O(A^2)$ as $A\to\infty$.
For ease of notation, we write
\beq\elab{eqn3.2}
A=\eps^{-1}~~\text{and}~~ B=\bar{B}\eps^{-2}
\eeq
with $\bar{B}=O(1)^+$ as $\eps\to 0$. 
In what follows, to avoid issues with the discontinuity in the cut-off nonlinearity, it proves convenient to make use of an exact result, obtained via \eref{KPP}, \eref{.noflux} and regularity on $u$, and given by,
\beq\elab{eqn3.3}
\begin{aligned}
\left(\int_0^1{u(x,s,t)}ds\right)_t + \eps^{-1}\left(\int_0^1{\alpha(s)u(x,s,t)}ds\right)_x = \left(\int_0^1{u(x,s,t)}ds\right)_{xx} + \int_0^1{f_c(u(x,s,t))}ds \\
 ~~\forall~~(x,t)\in \mathbb{R}\times\mathbb{R}^+.
 \end{aligned}
 \eeq
Moreover, since the cut-off nonlinearity {does not appear, at least up to order $O(\epsilon)$,} 
we are able to work, to our advantage, {directly} with {(IBVP)}, rather than the moving boundary formulation. We now introduce the expansion,
\beq\elab{eqn3.4}
u(x,y,t;\eps) = u_0(x,y,t) + \eps u_1(x,y,t) + O(\eps^2)
\eeq
as $\eps\to 0$, with $(x,y,t)\in \mathbb{R}\times[0,1]\times[0,\infty)$. On substitution from \eref{eqn3.4} into (IBVP), we obtain, at leading order, the following problem for $u_0(x,y,t)$, namely,
\begin{subequations}
\beq\elab{eqn3.5}
u_{0yy} = 0,~~(x,y,t)\in \mathbb{R}\times(0,1)\times(0,\infty),
\eeq
which must be solved subject to
\beq\elab{eqn3.6}
u_{0y}=0,~~(x,y,t)\in \mathbb{R}\times\{0,1\}\times(0,\infty),
\eeq
\end{subequations}
which requires,
\beq\elab{eqn3.7}
u_0(x,y,t) = \bar{u}(x,t),~~(x,y,t)\in \mathbb{R}\times[0,1]\times[0,\infty),
\eeq
with $\bar{u}:\mathbb{R}\times[0,\infty)\to \mathbb{R}$ to be determined subject to suitable regularity, with
\beq\elab{eqn3.8}
\bar{u}(x,t)\to
\begin{cases}
1~~\text{as}~~x\to -\infty, \\
0~~\text{as}~~x\to \infty,
\end{cases}
\eeq
uniformly for $t\in [0,T]$, any $T>0$, and
\beq\elab{eqn3.9}
\bar{u}(x,0) = H(-x),~~x\in \mathbb{R}.
\eeq
At $O(\eps)$ in (IBVP) we obtain from \eref{KPP}
\begin{subequations}\elab{eqn3.10-3.11}
\beq\elab{eqn3.10}
u_{1yy} = \frac{1}{\bar{B}}\alpha(y)\bar{u}_x(x,t),~~~~(x,y,t)\in \mathbb{R}\times(0,1)\times(0,\infty),
\eeq
which must be solved subject to
\beq\elab{eqn3.11}
u_{1y}=0,~~(x,y,t)\in \mathbb{R}\times\{0,1\}\times(0,\infty).
\eeq
\end{subequations}
The solution to {\eref{eqn3.10-3.11}} is given by,
\beq\elab{eqn3.12}
u_1(x,y,t) = \bar{V}(x,t) +\frac{1}{\bar{B}}\widetilde{\alpha}(y)\bar{u}_x(x,t)~~~~(x,y,t)\in \mathbb{R}\times[0,1]\times(0,\infty).
\eeq
Here $\bar{V}:\mathbb{R}\times[0,\infty)\to \mathbb{R}$ is to be determined, with suitable regularity, and
\beq\elab{eqn3.13}
\widetilde{\alpha}(y) = \int_0^y\left(\int_0^{\mu}{\alpha(s)}ds\right)d\mu~~\forall~~y\in [0,1].
\eeq
We observe from \eref{eqn3.13} that
\beq\elab{eqn3.14}
\widetilde{\alpha}'(y) = \int_0^y{\alpha(s)}ds,~~y\in [0,1],
\eeq
so that,
$\widetilde{\alpha}'(0)=\widetilde{\alpha}'(1)=0$ 
via \eref{eqn1.4}. We now substitute from \eref{eqn3.4}, \eref{eqn3.7}, \eref{eqn3.12} and \eref{eqn3.13} into \eref{eqn3.3}, and perform the integrations to obtain, at {$O(1)$}, the following PDE for $\bar{u}(x,t)$, namely,
\begin{subequations}\elab{eqn3.163.17}
\beq\elab{eqn3.16}
\bar{u}_t = (1+\bar{B}^{-1}\Delta)\bar{u}_{xx} + f_c(\bar{u}),~~(x,t)\in \mathbb{R}\times(0,\infty),
\eeq
with
\beq\elab{eqn3.17}
\Delta = \int_0^1{(\widetilde{\alpha}'(s))^2}ds > 0
\eeq
\end{subequations}
when $\alpha(y)$ is nontrivial. Equation \eref{eqn3.16} must be solved subject to the initial and boundary conditions \eref{eqn3.9} and \eref{eqn3.8}. With the scaling transformation
{$X = (1+\bar{B}^{-1}\Delta)^{-\frac{1}{2}}x$}. 
the initial value problem \eref{eqn3.16}, \eref{eqn3.9} and \eref{eqn3.8} for $\bar{u}$ reduces precisely to that studied in detail recently in \cite{Tisbury_etal2020a,Tisbury_etal2020b}. In \cite{Tisbury_etal2020a} it is established that the PDE \eref{eqn3.163.17} has a unique PTW solution with propagation speed
\beq\elab{eqn3.19}
v = \sqrt{(1 + \bar{B}^{-1}\Delta)}v^{\ast}(u_c),
\eeq
with $v^{\ast}(u_c)$ as given in \cite{Tisbury_etal2020a}.
It should be noted that $1 + \bar{B}^{-1}\Delta=\kappa_{\text{eff}}$, the effective diffusivity introduced in subsection 1.1.
 The initial boundary value problem \eref{eqn3.16}, \eref{eqn3.9} and \eref{eqn3.8} is studied in detail in \cite{Tisbury_etal2020b}, where it is principally established that $\bar{u}$ evolves into a PTW structure as $t\to \infty$.
We are now able to interpret this for (IBVP) via \eref{eqn3.4}. At \emph{leading order}, as $A\to \infty$ with $\bar{B} = (B/A^2)=O(1)$, we observe that $u$ becomes rapidly homogeneous in $y$ (on a time scale $t=O(A^{-2})$ as $A\to \infty$) and thereafter is governed by the one-dimensional cut-off KPP reaction-diffusion equation \eref{eqn3.16}. The effect of the advective shear flow is simply to enhance the streamwise diffusion coefficient from unity to $(1 + (A^2/B)\Delta)$. For each $A$, $B$, a unique PTW exists, which has propagation speed,
\beq\elab{eqn3.20}
\hat{v}(A,B,u_c) = (1+\bar{B}^{-1}(A,B)\Delta)^{\frac{1}{2}}v^{\ast}(u_c) + o(1)
\eeq
as $A\to \infty$, with $\bar{B} = B/A^2 = O(1)$ and $v^{\ast}(u_c)$ as given in \cite{Tisbury_etal2020a}. This PTW solution forms the large-$t$ attractor for (IBVP).

{As $B$ reduces in order relative to $A=\epsilon^{-1}$, an examination of the PDE \eref{KPP} reveals that the balancing of terms at $O(\epsilon)$ undergoes a change, in particular when the distinguished limit becomes $B=O(\eps^{-\frac{2}{n}})$ with $n>1$} . In this case we write
\beq\elab{eqn3.21}
B = \bar{B}\eps^{-\frac{2}{n}}
\eeq
with, again, $\bar{B}=O(1)$ as $\eps\to 0$. Returning to PDE \eref{KPP}, the natural scale in the $x$ coordinate is now of $O(\eps^{-(1-\frac{1}{n})})$ (so that the front thickness is increasing with $n$). Thus we introduce the coordinate $X$, with 
\beq\elab{eqn3.22}
x = \bar{B}^{\frac{1}{2}(n-1)}\eps^{-(1-\frac{1}{n})}X
\eeq
and the $O(1)$ factor is for algebraic convenience. We next follow the earlier case, when $n=1$, and expand as,
\beq\elab{eqn3.23}
u(X,y,t;\eps) = \bar{u}(X,t) + O(\eps^{\frac{1}{n}})
 \eeq
as $\eps\to 0$ with $(X,y,t)\in \mathbb{R}\times[0,1]\times[0,\infty)$. Here $\bar{u}:\mathbb{R}\times[0,\infty)\to \mathbb{R}$ now satisfies,
\beq\elab{3.24}
\bar{u}_t = \bar{B}^{-n}\Delta\bar{u}_{XX} + f_c(\bar{u}),~~(X,t)\in \mathbb{R}\times(0,\infty),
\eeq
with suitable regularity, and again subject to initial and boundary conditions \eref{eqn3.9} and \eref{eqn3.8}.The conclusions regarding (IBVP) and (BVP) are thus as for the case $n=1$, \emph{except} now the front is on the stretched length scale $x=O(A^{(1-\frac{1}{n})})$ and \eref{eqn3.20} is modified to,
\beq\elab{eqn3.25}
\hat{v}(A,B,u_c) = \frac{\sqrt{\Delta}}{\sqrt{\bar{B}(A,B)}}A^{(1-\frac{1}{n})}v^{\ast}(u_c) + o(A^{(1-\frac{1}{n})})
\eeq
as $A\to \infty$, with $\bar{B}(A,B)=B/A^{\frac{2}{n}}=O(1)$. 
This completes the structure to (IBVP) and (BVP) in the case of slowly varying fronts, with $A\gg1$ and $B\gg1$. 

We end this section with  comments on the behaviour of the 
PTW propagation speed. 
We observe from 
\eref{eqn3.20} and \eref{eqn3.25} 
that $\hat{v}(A,B,u_c)$   is enhanced by the shear flow through a prefactor that is entirely determined by the flow whilst the effects of reaction cut-off are felt through the factor $v^*(u_c)$ corresponding to 
 the propagation speed for the reaction cut-off problem in the absence of a flow. Thus, the asymptotic limits \eref{asympt1} and \eref{asympt2} are, in the presence of a shear flow, magnified by a factor entirely dependent on the flow. 
 For $B=O(A^{\frac{2}{n}})$, the enhancement is significant. 
Finally, simple calculations for the Couette and Poiseulle flows give, specifically,
\beq\elab{Delta_example}
\Delta=\begin{cases}
    1/120, &\quad \text{for}\quad \alpha(y)=y-1/2~(\text{Couette}),\\
    1/1890, &\quad \text{for}\quad \alpha(y)=-2y^2+2y-1/3 ~(\text{Poiseulle}). 
\end{cases}
\eeq 
Thus, the PTW propagates faster in the plane Couette flow than in the plane Poiseulle flow.

\section{Slowly varying, balanced or rapidly varying front with weak advection: $B^{-1}=o(1), B=O(1)$ or $B=o(1)$ as $A\to 0$}
In this section we consider (IBVP) ((QIVP)) and (BVP ((QBVP)) in the case when the advection velocity scale is small compared to advectionless front propagation speed, so that $A\ll1$. To begin with we examine the structure to (IBVP) when also $B\gg1$, so that the channel width is small compared to the advectionless front thickness. To formalise this case, we consider (IBVP) when $B^{-1}=o(1)$ as $A\to 0$.
\subsection{Slowly varying front: $B^{-1}=o(1)$ as $A\to 0$}
This limit considers the situation when the channel width is small compared to the thickness of the advectionless front. For notational convenience, we write
\beq\elab{eqn4.1}
A = \delta
\eeq
and then {a balance in \eref{eqn2.10} determines an expansion}  in the form,
\beq\elab{eqn4.2}
u(x,y,t;\delta) = u_0(x,y,t) + O(B^{-1}(\delta))
\eeq
as $\delta\to 0$ with $(x,y,t)\in \mathbb{R}\times[0,1]\times[0,\infty)$. {Here, formally, $B^{-1}(\delta)=o(1)$ as $\delta\to 0$.} On substituting from \eref{eqn4.2} into (IBVP), we obtain,
\beq\elab{eqn4.3}
u_0(x,y,t) = \bar{u}(x,t),~~(x,y,t)\in \mathbb{R}\times[0,1]\times[0,\infty),
\eeq
with $\bar{u}:\mathbb{R}\times[0,\infty)\to \mathbb{R}$ to be determined, with suitable regularity. Moving to $O(B^{-1}(\delta))$, we proceed as in section 3, and, without repeating details, we obtain the following scalar problem for $\bar{u}(x,t)$, namely,
\beq\elab{eqn4.4}
\bar{u}_t = \bar{u}_{xx} + f_c(\bar{u}),~~(x,t)\in \mathbb{R}\times(0,\infty),
\eeq
together with initial and boundary conditions \eref{eqn3.9} and \eref{eqn3.8}. This is precisely the scalar evolution problem studied in detail in \cite{Tisbury_etal2020a,Tisbury_etal2020b}, where it is established that $\bar{u}$ evolves into a PTW structure as $t\to \infty$, and, for each $u_c\in (0,1)$, this PTW solution is unique, and has propagation speed $v^{\ast}(u_c)$, as given in \cite{Tisbury_etal2020a}. Thus in this case, for (BVP),
\beq\elab{eqn4.5}
\hat{v}(A,B,u_c) = v^{\ast}(u_c) + O\left(A^2B^{-1}(A)\right)
\eeq
as $A\to 0$ with $B^{-1}(A)=o(1)$. This simple case is now complete. To summarise this case, it follows from \eref{eqn4.3} that the deformation of the front interface is weak 
and of order $O(AB^{-1}(A))$ whilst from {\eref{eqn4.5}} the correction to the propagation speed is even weaker. This correction can be  considered via restricting attention to the problem (QBVP) since at higher order the detailed nature of the discontinuity in the reaction must be addressed and so, the free-boundary formulation (QBVP) needs to be adopted. We now move on to the next case.

\subsection{Balanced front: $B=O(1)$ as $A\to 0$}\slab{balanced}
This limit addresses the situation when the channel width is comparable to the  front thickness in the absence of advection. In considering the limit $B=O(1)$ as $A\to 0$, since the discontinuity in the cut-off nonlinearity is encountered immediately, at leading order, it is most convenient to address the formulation (QIVP), and we restrict attention to PTW solutions to (QIVP), via the formulation (QBVP). We look for a solution to (QBVP) expanded in the form (with $A=\delta$),
\begin{subequations}\elab{eqn4.64.74.8}
\beq\elab{eqn4.6}
U_T(\xi',y';\delta) = u_0(\xi',y') + \delta u_1(\xi',y') + O(\delta^2)
\eeq
as $\delta\to 0$, with $(\xi',y')\in Q_L\cap Q_R$. In addition, for $y'\in [0,L]$, we write,
\beq\elab{eqn4.7}
\zeta(y',\delta) = \zeta_0(y') + \delta\zeta_1(y') + O(\delta^2)
\eeq
and expand the propagation speed,
\beq\elab{eqn4.8}
v(\delta) = v_0 + \delta v_1 + O(\delta^2).
\eeq
\end{subequations}\elab{eqn4.9-4.15}
We now substitute from \eref{eqn4.64.74.8} into (QBVP). Collecting terms at $O(1)$, we obtain the following problem for $u_0$, $\zeta_0$ and $v_0$, namely,
\begin{subequations}\elab{eqn4.94.15}
\beq\elab{eqn4.9}
\nabla^2u_0 + v_0u_{0\xi'} + f_c(u_0) = 0,~~(\xi',y')\in Q_L\cup Q_R
\eeq
with
\beq\elab{eqn4.10}
u_0\ge u_c~\text{in}~Q_L,~~u_0\le u_c~\text{in}~Q_R,
\eeq
and
\beq\elab{eqn4.11}
u_0(\xi',y')\to
\begin{cases}
1~~\text{as}~~\xi'\to -\infty,\\
0~~\text{as}~~\xi'\to \infty,
\end{cases}
\eeq
uniformly for $y'\in [0,L]$, whilst
\beq\elab{eqn4.12}
u_{0y'}(\xi',0)=u_{0y'}(\xi',L)=0,~~\xi'\in \mathbb{R}.
\eeq
\beq\elab{eqn4.13}
u_0(\zeta_0(y')^+,y')=u_0(\zeta_0(y')^-,y')=u_c,~~y'\in [0,L],
\eeq
\beq\elab{eqn4.14}
[(u_{0\xi'}-\zeta_{0y'}u_{0y'})(\zeta_0(y'),y')]_L^R=0,~~y'\in (0,L).
\eeq
\beq\elab{eqn4.15}
\int_0^L{\zeta_0(y')}dy'=0.
\eeq
\end{subequations}
The elliptic problem {\eref{eqn4.94.15}} has a unique solution, and this solution is independent of $y'$. The solution is given in \cite{Tisbury_etal2020a} (see Theorem 1) in terms of $U_T:\mathbb{R}\to \mathbb{R}$, namely,
\begin{subequations}\elab{eqn4.164.18}
\beq\elab{eqn4.16}
u_0(\xi',y') = U_T(\xi'),~~(\xi',y')\in \mathbb{R}\times [0,L]
\eeq
with
\beq\elab{eqn4.17}
\zeta_0(y')=0~~\forall~~y'\in [0,L]
\eeq
and
\beq\elab{eqn4.18}
v_0=v^{\ast}(u_c).
\eeq
\end{subequations}
We recall, from \cite{Tisbury_etal2020a}, that $U_T(\xi')$ is {monotone} decreasing in $\xi'$, with,
\begin{subequations}
\beq\elab{4.19}
U_T''(0^+)-U_T''(0^-)=-f_c,
\eeq
\beq\elab{4.20}
U_T(\xi') = u_ce^{-v^{\ast}(u_c)\xi'}~~\forall~~\xi'\in [0,\infty),
\eeq
\beq\elab{4.21}
U_T(\xi') \sim 1 - A_{-\infty}e^{\lambda_+(v^{\ast}(u_c))\xi'}~~\text{as}~~\xi'\to -\infty,
\eeq
with the constant $A_{-\infty}>0$ and
\beq\elab{4.22}
\lambda_+(v) = \frac{1}{2}\left((v^2 + 4|f_c'(1)|)^{\frac{1}{2}} - v\right).
\eeq
\end{subequations}
We now proceed to terms at $O(\delta)$ in (QBVP). This results in the following inhomogeneous linear elliptic boundary value problem for $u_1$, $\zeta_1$ and $v_1$, namely,
\begin{subequations}\elab{eqn4.234.27}
\beq\elab{eqn4.23}
\nabla^2u_1 + v^{\ast}u_{1\xi'} + H(-\xi')f_c'(U_T(\xi'))u_1 = (\bar{\alpha}(y')-v_1)U_T'(\xi'),~~(\xi',y')\in (\mathbb{R}\setminus \{0\})\times (0,L),
\eeq
\beq\elab{eqn4.24}
u_1(\xi',y')\to 0~~\text{as}~~|\xi'|\to \infty~~\text{uniformly for}~~y\in [0,L],
\eeq
\beq\elab{eqn4.25}
u_{1y'}(\xi',0)=u_{1y'}(\xi',L)=0,~~\xi'\in \mathbb{R},
\eeq
\beq\elab{eqn4.26}
u_1(0^+,y')=u_1(0^-,y')=v^{\ast}u_c\zeta_1(y'),~~y'\in (0,L),
\eeq
\beq\elab{eqn4.27}
u_{1\xi'}(0^+,y') - u_{1\xi'}(0^-,y') = -f_c\zeta_1(y'),~~y'\in (0,L),
\eeq
where \eref{eqn4.26} and \eref{eqn4.27}
were obtained using \eref{4.19} and \eref{4.20}.
\end{subequations}
We first eliminate $\zeta_1(y')$ from the problem \eref{eqn4.24}-\eref{eqn4.27}, which reduces the boundary conditions \eref{eqn4.26} and \eref{eqn4.27} to
\begin{subequations}\elab{eqn4.284.29}
\beq\elab{eqn4.28}
u_1(0^+,y')=u_1(0^-,y'),
\eeq
\beq\elab{eqn4.29}
u_{1\xi'}(0^+,y') - u_{1\xi'}(0^-,y') = -f_c(v^{\ast}u_c)^{-1}u_1(0^-,y'),
\eeq
\end{subequations}
for all $y'\in (0,L)$, after which $\zeta_1(y')$ is given by
\beq\elab{eqn4.30}
\zeta_1(y') = (v^{\ast} u_c)^{-1}u_1(0^-,y'),~~y'\in [0,L].
 \eeq
The linear elliptic problem {\eref{eqn4.234.27} and \eref{eqn4.284.29}}  can be solved explicitly. We begin by using Fourier's Theorem to write,
\beq\elab{eqn4.31}
u_1(\xi',y') = \sum_{n=0}^{\infty}a_n(\xi')\text{cos}\left(\frac{n\pi}{L}y'\right),~~(\xi',y')\in \mathbb{R}\times[0,L],
\eeq
with $a_n:\mathbb{R}\to \mathbb{R}$ ($n=0,1,2,...$) to be determined. We observe that the boundary conditions \eref{eqn4.25} are satisfied by \eref{eqn4.31}. In addition, we may write, via \eref{eqn2.21},
\begin{subequations}\elab{eqn4.314.32}
\beq\elab{eqn4.31'}
\bar{\alpha}(y') = \sum_{n=1}^{\infty}\bar{\alpha}_n\text{cos}\left(\frac{n\pi}{L}y'\right),~~y'\in [0,L],
\eeq
with
\beq\elab{eqn4.32}
\bar{\alpha}_n = \frac{2}{L}\int_0^L{\bar{\alpha}(s)\text{cos}\left(\frac{n\pi}{L}s\right)}ds
\eeq
\end{subequations}
for $n=1,2,\ldots$. Now substitute from \eref{eqn4.31} and \eref{eqn4.31'} into \eref{eqn4.23}--\eref{eqn4.26}, \eref{eqn4.28} and \eref{eqn4.29}. At $n=0$ we obtain the following problem for $a_0(\xi')$,
\begin{subequations}\elab{eqn4.334.36}
\beq\elab{eqn4.33}
a_0'' + v^{\ast}a_0' + H(-\xi')f'(U_T(\xi'))a_0 = -v_1U_T'(\xi'),~~\xi'\in \mathbb{R} \setminus \{0\},
\eeq
\beq\elab{eqn4.34}
a_0(0^+)=a_0(0^-),
\eeq
\beq\elab{eqn4.35}
a_0'(0^+) - a_0'(0^-) = -f_c(v^{\ast}u_c)^{-1}a_0(0^-),
\eeq
\beq\elab{eqn4.36}
a_0(\xi')\to 0~~\text{as}~~|\xi'|\to \infty.
\eeq
\end{subequations}
To solve {\eref{eqn4.334.36}}, we first observe that $U_T'(\xi')$ is a solution to the homogeneous form of the linear ODE \eref{eqn4.33}. After writing this ODE in self-adjoint form, it is then readily established that the linear, inhomogeneous, boundary value problem \eref{eqn4.33}--\eref{eqn4.36} has the \emph{solvability condition}
\beq\elab{eqn4.36'}
v_1\left(\int_{-\infty}^0{(U_T'(s))^2e^{v^{\ast}s}}ds + v^{\ast}u_c^2 \right) = 0
\eeq
which requires
\beq\elab{eqn4.37}
v_1 = 0.
\eeq
The solution to {\eref{eqn4.334.36}} is then given by,
\beq\elab{eqn4.38}
a_0(\xi') = CU_T'(\xi'),~~\xi'\in \mathbb{R},
\eeq
with $C$ a constant to be determined. Now, from \eref{eqn2.67} and \eref{eqn4.7}, we require
\beq\elab{eqn4.39}
\int_0^L{\zeta_1(y')}dy' = 0,
\eeq
which becomes, via \eref{eqn4.30},
\beq\elab{eqn4.39'}
\int_0^L{u_1(0^-,y')}dy' = 0
\eeq
and so, via \eref{eqn4.31} and \eref{eqn4.34},
\beq\elab{eqn4.40}
a_0(0^+)=a_0(0^-)=0.
\eeq
Now $U_T'(0) = -v^{\ast}u_c<0$, and so \eref{eqn4.38} and \eref{eqn4.40} requires $C=0$, and therefore
\beq\elab{eqn4.41}
a_0(\xi') = 0,~~\xi'\in \mathbb{R}.
\eeq
For $n=1,2,3...$, we obtain the problem,
\begin{subequations}\elab{eqn4.424.45}
\beq\elab{eqn4.42}
a_n'' + v^{\ast}a_n' - (n^2\pi^{2}L^{-2} -  H(-\xi')f'(U_T(\xi')))a_n = \bar{\alpha}_nU_T'(\xi'),~~\xi'\in \mathbb{R} \setminus \{0\},
\eeq
\beq\elab{eqn4.43}
a_n(0^+)=a_n(0^-),
\eeq
\beq\elab{eqn4.44}
a_n'(0^+) - a_n'(0^-) = -f_c(v^{\ast}u_c)^{-1}a_n(0^-),
\eeq
\beq\elab{eqn4.45}
a_n(\xi')\to 0~~\text{as}~~|\xi'|\to \infty.
\eeq
\end{subequations}
The solution to {\eref{eqn4.424.45}} is readily obtained as,
\beq\elab{eqn4.46}
a_n(\xi') = -\bar{\alpha}_nL^2(n^2\pi^2)^{-1}U_T'(\xi'),~~\xi'\in \mathbb{R}.
\eeq
Thus, via \eref{eqn4.41} and \eref{eqn4.46}, we have, from \eref{eqn4.31},
\beq\elab{eqn4.47}
u_1(\xi',y') = - \frac{U_T'(\xi')L^2}{\pi^2}\sum_{n=1}^{\infty}\frac{\bar{\alpha}_n}{n^2}\text{cos}\left(\frac{n\pi}{L}y'\right),~~(\xi',y')\in \mathbb{R}\times[0,L].
\eeq
 We recall that $\bar{\alpha}(y')$ is given by {\eref{eqn4.314.32}}, and so,
it is convenient to introduce $\bar{a}(y')$ as,
\beq\elab{eqn4.47'}
\bar{a}(y') \equiv \int_0^{y'} {\bar{\alpha}(s)}ds \\
=\frac{L}{\pi} \sum_{n=1}^{\infty}\frac{\bar{\alpha}_n }{n}\text{sin}\left(\frac{n\pi}{L}y'\right),~~y'\in [0,L].
 \eeq
{Integrating \eref{eqn4.47'} and rearranging yields} 
\beq\elab{eqn4.48}
\frac{L^2}{\pi^2}
\sum_{n=1}^{\infty}\frac{\bar{\alpha}_n }{n^2}\text{cos}\left(\frac{n\pi}{L}y'\right) = \frac{L^2}{\pi^2}\sum_{n=1}^{\infty}\frac{\bar{\alpha}_n}{n^2} - \int_0^{y'} {\bar{a}(s)}ds  ,~~y'\in [0,L].
\eeq
Thus we have,
\beq\elab{eqn4.49}
u_1(\xi',y') = \left(\int_0^{y'}{\bar{a}(s)}ds - \frac{L^2}{\pi^2}\sum_{n=1}^{\infty}{\frac{\bar{\alpha}_n}{n^2}}\right)U_T'(\xi'),~~(\xi',y')\in \mathbb{R}\times[0,L].
\eeq
Now, from \eref{eqn4.48}, we obtain,
\beq\elab{eqn4.50}
\frac{1}{L}\int_0^L \left(\int_0^{y'}{\bar{a}(s)}ds\right)dy' = \frac{L^2}{\pi^2}\sum_{n=1}^{\infty}{\frac{\bar{\alpha}_n}{n^2}}.
\eeq
It is thus convenient to introduce
\beq\elab{eqn4.51}
\phi(y') = \int_0^{y'}{\bar{a}(s)}ds = \int_0^{y'}\left(\int_0^{w} {\bar{\alpha}(s)}ds\right)dw,~~y'\in [0,L].
\eeq
We can now use \eref{eqn4.50} and \eref{eqn4.51} to re-write \eref{eqn4.49} as,
\beq\elab{eqn4.52}
u_1(\xi',y') = \left(\phi(y') - \langle\phi\rangle_L\right)U_T'(\xi'),~~(\xi',y')\in \mathbb{R}\times[0,L],
\eeq
with $\langle\cdot\rangle_L$ denoting the usual mean value on the interval $[0,L]$. Finally, from \eref{eqn4.30} and \eref{eqn4.52}, we obtain,
\beq\elab{eqn4.53}
\zeta_1(y') = \langle\phi\rangle_L - \phi(y'),~~y'\in [0,L].
\eeq
Thus, in this case we have, via 
{\eref{eqn4.64.74.8},
\eref{eqn4.164.18}},  
\eref{eqn4.37}, \eref{eqn4.41}, \eref{eqn4.52} and \eref{eqn4.53},
\begin{subequations}
\beq\elab{eqn4.54}
U_T(\xi',y',\delta) = U_T(\xi') + \delta\left(\phi(y') - \langle\phi\rangle_L\right)U_T'(\xi') + O(\delta^2),~~(\xi',y')\in \mathbb{R}\times[0,L],
\eeq
\beq\elab{eqn4.55}
\zeta(y',\delta) = \delta(\langle\phi\rangle_L - \phi(y')),~~y'\in [0,L].
\eeq
\beq\elab{eqn4.56}
v(\delta) = v^{\ast}(u_c) + O(\delta^2),
\eeq
\end{subequations}
as $\delta\to 0$ with $L=B^{-\frac{1}{2}}=O(1)^+$. We conclude in this case that there is a unique PTW solution, given by \eref{eqn4.54} and \eref{eqn4.55}, for each $u_c\in (0,1)$, which has propagation speed,
\beq\elab{eqn4.57}
\hat{v}(A,B,u_c) = v^{\ast}(u_c) + O(A^2)
\eeq
as $A\to 0$ with $B=O(1)$. We note that, although the correction to the PTW propagation speed is at least as small as $O(A^2)$, the corrections to the PTW structure are both at $O(A)$.
We also note from \eref{eqn4.55} with \eref{eqn4.51} and \eref{eqn4.47'} that $\zeta(y',\delta)$ has interior stationary points if and only if $\bar a(y')$ has a zero in the interior of the domain.
 This completes the asymptotic analysis in the case $A\to 0$ with $B=O(1)$. 
Contrasting with the previous case, the correction to the propagation speed is of the same order when $B(A)=O(1)$ as $A\to 0$, which is as expected. In fact we see that in both cases the correction to the propagation speed can be written, uniformly as $O(A^2/(B(A)+A))$ as $A\to 0$, and we will see that this continues to hold in each of the cases considered below in the rest of this section.
We finally compute the interface deformation for the Couette flow \eref{Couette} and the Poiseuille flow \eref{Poiseuille} using \eref{eqn4.55} and \eref{eqn4.51}.  
Figure \fref{interface_smallA} compares  the deformation against the structure of the two flows.
It is clear that the deformation  in the Couette flow is more pronounced than in the plane Poiseulle flow.

\begin{figure}
		 \begin{center}
 	     \includegraphics[width=1\textwidth]{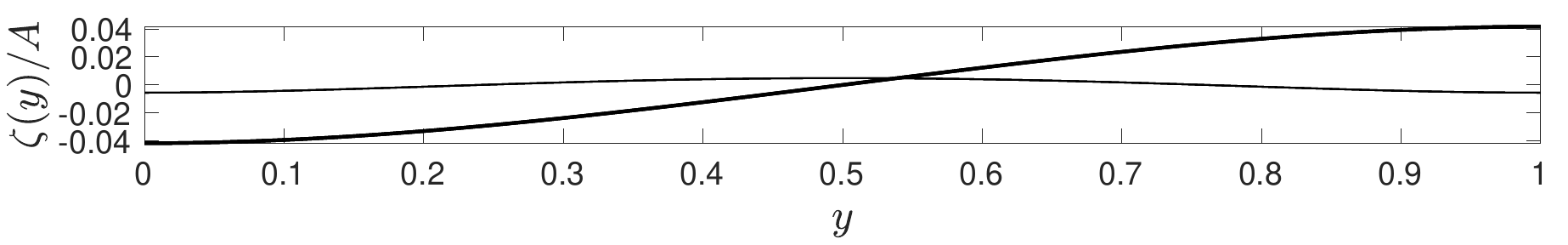}
		  \end{center}
    \caption {A graph of  the interface $\zeta(y)/A$ as $A\to 0$ obtained  for $B=1$ 
	using \eref{eqn4.55} and \eref{eqn4.51} 
	for
	     the Couette     (thick solid lines) and the 
	     Poiseuille   flow (thin solid lines). 
}
  \flab{interface_smallA}  
  \end{figure}

The final case in this section is $A\ll1$ with $B\ll1$, and the appropriate distinguished limit is $B=O(A)$ as $A\to 0$.

\subsection{Rapidly varying front: $B=O(A)$ as $A\to 0$}
Here we consider the asymptotic structure of PTW solutions when the front thickness in the absence of advection is small compared to the channel width. In this case, it is natural to consider PTW solutions via the formulation (BVP). Again, setting $A=\delta$, we write,
\beq\elab{eqn4.58}
B = \bar{B}\delta
\eeq
with $\bar{B}=O(1)$ as $\delta\to 0$. We now expand in the form,
\begin{subequations}\elab{eqn4.594.61}
\beq\elab{eqn4.59}
U_T(z,y;\delta) = u_0(z,y) + \delta u_1(z,y) + O(\delta^2),~~(z,y)\in \mathbb{R}\times[0,1],
\eeq
\beq\elab{eqn4.60}
\gamma(y;\delta) = \gamma_0(y) + \delta\gamma_1(y) + O(\delta^2),~~y\in [0,1],
\eeq
\beq\elab{eqn4.61}
v(\delta) = v_0 + \delta v_1 + O(\delta^2)
\eeq
\end{subequations}
as $\delta\to 0$. We substitute \eref{eqn4.58}  and \eref{eqn4.594.61} into (BVP).
At leading order we obtain the following problem for $u_0$, $v_0$ and $\gamma_0$, namely,
\begin{subequations}\elab{eqn4.624.67}
\beq\elab{eqn4.62}
u_{0zz} + v_0u_{0z} + f_c(u_0) = 0,~~(z,y)\in \mathbb{R}\times(0,1),
\eeq
\beq\elab{eqn4.63}
u_0(\gamma_0^{\pm}(y),y)=u_c,~~y\in [0,1],
\eeq
\beq\elab{4.64}
u_{0z}(\gamma_0^+(y),y) = u_{0z}(\gamma_0^-(y),y),~~y\in [0,1],
\eeq
\beq\elab{eqn4.65}
u_0(z,y)\to
\begin{cases}
1~~\text{as}~~z\to -\infty,\\
0~~\text{as}~~z\to \infty,
\end{cases}
\eeq
uniformly for $y\in [0,1]$,
\beq\elab{eqn4.66}
u_{0y}(z,0) = u_{0y}(z,1) = 0,~~z\in \mathbb{R},
\eeq
\beq\elab{eqn4.67}
\int_0^1{\gamma_0(y)}dy = 0.
\eeq
\end{subequations}
{Since partial derivatives with respect to $y$ do not appear in the problem  \eref{eqn4.62}--\eref{eqn4.65} we can address this problem following 
 \cite{Tisbury_etal2020a} and regard $y$ as a parameter. Thus, we may conclude that, up to translation invariance,  \eref{eqn4.624.67} has a unique solution. Since we may regard the translational invariance as dependent on $y$ then we may write the solution  as} 
\beq\elab{eqn4.68}
u_0(z,y) = U_T(z-\gamma_0(y)),~~(z,y)\in \mathbb{R}\times[0,1],
\eeq
\beq\elab{eqn4.69}
v_0=v^{\ast}(u_c),
\eeq
with $\gamma_0:[0,1]\to \mathbb{R}$ to be determined, so that $\gamma_0\in C^1([0,1])$, and satisfies  condition \eref{eqn4.67}, together with,
\beq\elab{eqn4.70}
\gamma_0'(0)=\gamma_0'(1)=0
\eeq
via \eref{eqn4.68} with \eref{eqn4.66}. We next formulate the problem at $O(\delta)$ from (BVP), which provides an inhomogeneous linear problem for $u_1$ and $\gamma_1$. After a significant amount of detailed, but routine, calculation on this problem (which, for brevity, we do not include here), we find that it requires a solvability condition {arising from the classical Fredholm alternative which must be satisfied for the linear inhomogenous boundary value problem to have a solution.} This provides a nonlinear eigenvalue problem which determines $\gamma_0:[0,1]\to \mathbb{R}$ and $v_1\in \mathbb{R}$, namely,
\beq\elab{eqn4.71}
\gamma_0'' + \frac{1}{2}v^{\ast}(\gamma_0')^2 + \bar{B}^{-1}(\alpha(y)-v_1) = 0,~~y\in (0,1),
\eeq
together with boundary conditions \eref{eqn4.70} and condition \eref{eqn4.67}, which we henceforth refer to as (EP). In considering (EP), we first introduce $\psi: [0,1]\to \mathbb{R}$ given by
\beq\elab{eqn4.72}
\psi(y) = \text{exp}\left(\frac{1}{2}v^{\ast}\gamma_0(y)\right)~~\forall~~y\in [0,1]. 
\eeq
{which requires that}
\beq\elab{eqn4.75}
\psi(y)>0,~~y\in [0,1]
\eeq
and, upon using \eref{eqn4.67}, that 
\beq\elab{eqn4.76}
\int_0^1{\log(\psi(y))}dy = 0. 
\eeq
In terms of this new dependent variable, (EP) 
{becomes a classical regular Sturm-Liouville eigenvalue problem (see, for example, Coddington and Levinson \cite{CoddingtonLevinson}) that we will henceforth refer to as (SL):} 
 \begin{taggedsubequations}{SL}\label{SL}
\beq\elab{eqn4.73}
\psi'' + (k\alpha(y) + \lambda)\psi = 0,~~y\in (0,1),
\eeq
{with boundary  conditions} 
 \beq\elab{eqn4.74}
 \psi'(0)=\psi'(1)=0,
 \eeq
 {and   condition} 
where 
\[
\lambda=-\frac{1}{2}v^{\ast}\bar{B}^{-1}v_1,~~k=\frac{1}{2}v^{\ast}\bar{B}^{-1}.
\]
\end{taggedsubequations}
  Indeed, {(SL)} has arisen via a very different route in the work of Haynes and Vanneste \cite{HaynesVanneste2014a} in the context of a purely advection-diffusion problem at high P\'eclet number (corresponding to high $A/\sqrt{B}$), and in relation to the present context, it has been studied by them for the particular Couette and plane Poiseuille flows noted in \eref{CouettePoiseuille}. The positivity requirement \eref{eqn4.75} 
  dictates  that we require the \emph{smallest (principal)} eigenvalue of this Sturm-Liouville problem.
  We denote the principal eigenvalue by
\beq\elab{eqn4.78}
\lambda = \lambda_0(k)
\eeq
and the associated principal, $L^1$-normalised, eigenfunction by $\psi=\psi_0:[0,1]\to \mathbb{R}$ with 
\beq\elab{eqn4.79}
\psi_0 = \psi_0(y,k)>0~~\forall~~y\in [0,1]
\eeq
and
\beq\elab{eqn4.80}
\int_0^1{\psi_0(y,k)}dy = 1.
\eeq
with the choice of $L^1$ normalisation being convenient at a later stage. The solution to (EP) is then given by, on satisfying the final condition \eref{eqn4.76},
\beq\elab{eqn4.81}
\gamma_0(y) = \frac{2}{v^{\ast}}\left(\log(\psi_0(y,k)) - \langle\log(\psi_0(\cdot,k))\rangle \right)~~\forall~~y\in [0,1],
\eeq
with
\beq\elab{eqn4.82}
v_1 = -\frac{\lambda_0(k)}{k}.
\eeq
We now consider (SL) in the two cases when $k\gg1$ and $k\ll1$.

\subsubsection{(SL) when $k\ll1$}
We consider (SL) as $k\to 0$. An exposition of this type has been developed in \cite{HaynesVanneste2014a}, and for completeness we give a brief development in the present context. It follows from {(SL) with \eref{eqn4.79} and 
\eref{eqn4.80}} that 
\beq\elab{eqn4.83'}
\lambda_0(k)\to 0,
\quad\text{and}\quad
\psi_0(y,k)\to 1~~\text{uniformly for}~~y\in [0,1],
\eeq
as $k\to 0$. Thus we expand in the form,
\begin{subequations}\elab{eqn4.844.85}
\beq\elab{eqn4.84}
\psi_0(y,k) = 1 + k\bar{\psi}_0(y) + O(k^2),~~y\in [0,1],
\eeq
\beq\elab{eqn4.85}
\lambda_0(k) = k\bar{\lambda}_0 + O(k^2),
\eeq
\end{subequations}
as $k\to 0$. On substitution into (SL) using \eref{eqn4.80},  
we obtain the problem
\begin{subequations}\elab{eqn4.864.88}
\beq\elab{eqn4.86}
\bar{\psi}_0'' + (\alpha(y) + \bar{\lambda}_0) = 0,~~y\in (0,1),
\eeq
\beq\elab{eqn4.87}
\bar{\psi}_0'(0)=\bar{\psi}_0'(1)=0,
\eeq
\beq\elab{eqn4.88}
\int_0^1{\bar{\psi}_0(y)}dy=0.
\eeq
\end{subequations}
The solution to \eref{eqn4.864.88} is readily obtained as
\begin{subequations}\elab{eqn4.894.90}
\beq\elab{eqn4.89}
\bar{\psi}_0(y) =\left( \langle\widetilde{\alpha}(\cdot)\rangle -~ \widetilde{\alpha}(y)\right)~~\forall~~y\in [0,1],
\eeq
with $\widetilde{\alpha}:[0,1]\to \mathbb{R}$ given by \eref{eqn3.13}, and,
\beq\elab{eqn4.90}
\bar{\lambda}_0=0.
\eeq
\end{subequations}
Thus, via \eref{eqn4.81}, \eref{eqn4.82}, \eref{eqn4.844.85} and \eref{eqn4.894.90}, we obtain,
\beq\elab{eqn4.93}
v_1 = O(k),
\eeq
\beq\elab{eqn4.94}
\gamma_0(y) = \frac{2k}{v^{\ast}(u_c)}\left(\langle\widetilde{\alpha}(\cdot)\rangle -~ \widetilde{\alpha}(y)\right) + O(k^2)~~\forall~~y\in [0,1],
\eeq
as $k\to 0$. It should be noted that these results, \emph{for $\bar{B}\gg1$}, are in full accord with those of section \sref{balanced}, \emph{for {$B=O(1)$}}. 
Expression \eref{eqn4.93} can be made more precise by  performing higher-order corrections to $\lambda_0(k)$, the details of which are presented in \cite{HaynesVanneste2014a} and give
\beq \lambda_0(k)=-k^2\left\langle\left(\int_{0}^y\alpha(s)ds\right)^2\right\rangle+O(k^3),
\eeq
as $k\to 0$. 
Upon using \eref{eqn4.82} 
these corrections yield  
\beq
v_1=k \left\langle\left(\int_{0}^y\alpha(s)ds\right)^2\right\rangle+O(k^2),
\eeq
as $k\to 0$.

\subsubsection{(SL) when $k\gg1$}
We consider (SL) as $k\to \infty$. An analysis of this type has been given for the Couette and plane Poiseuille flow in \cite{HaynesVanneste2014a}. Here we develop the theory to apply to any shear flow which has a continuous derivative on $[0,1]$, and whose \emph{absolute maximum on $[0,1]$ occurs strictly at an interior point}. First, let
\beq\elab{eqn4.95}
\alpha_M = \sup_{y\in [0,1]}\alpha(y)>0
\eeq
when $\alpha(y)$ is nontrivial, via \eref{eqn1.4}. Now, we   restrict attention to the situation when $\alpha_M$ is achieved ($\alpha\in C^1([0,1])$) at a \emph{ single point, which is in the interior of $[0,1]$}. We label this point as $y_M\in (0,1)$. 
Numerical experiments on (SL) with $k\gg1$ suggest that the principal eigenfunction, $\psi_0$, under the normalisation \eref{eqn4.80}, 
develops into a $\delta$-sequence as $k\to \infty$, based at $y=y_M$. {This leads us to write} 
\beq\elab{eqn4.98}
\psi_0(y,k) \sim \delta(y-y_M),~~y\in [0,1]
\eeq
as $k\to \infty$, which automatically satisfies the {normalisation \eref{eqn4.80}}. Now, upon integration and using {\eref{eqn4.80}}, (SL) becomes 
\[
k\int_0^1{\alpha(y)\psi_0(y,k)}dy + \lambda_0(k) = 0
\]
which gives, using \eref{eqn4.98},
{$\lambda_0(k) \sim -\alpha_Mk$ as $k\to \infty$}.
This suggests the principal eigenvalue of (SL) should be expanded as,
\beq\elab{eqn4.100}
\lambda_0(k) = -\alpha_Mk + \tilde{\lambda}k^{\frac{1}{2}} + O(1)
\eeq
as $k\to \infty$, with $\widetilde{\lambda}\in \mathbb{R}$ to be determined. Using this, we may now determine the structure of $\psi_0:[0,1]\to \mathbb{R}$ as $k\to \infty$, in more detail, together with determining $\widetilde{\lambda}$. We first anticipate, from {(SL)}, 
that there will develop boundary layers at $y=0$ and $y=1$, as $k\to \infty$. We focus attention first at $y=0$. The boundary layer thickness will be $y=O(k^{-\frac{1}{2}})$, and so we introduce
\beq\elab{eqn4.101}
Y = k^{\frac{1}{2}}y = O(1)^+
\eeq
as $k\to \infty$ in the boundary layer. At leading order, (SL) becomes
\begin{subequations}\elab{eqn4.1024.102'}
\beq\elab{eqn4.102}
\psi_{0YY} - (\alpha_M -\alpha(0))\psi_0 = 0,~~Y>0,
\eeq
\beq\elab{eqn4.102'}
\psi_{0Y}(0)=0.
\eeq
\end{subequations}
The solution to {\eref{eqn4.1024.102'}} has,
\beq\elab{eqn4.103}
\psi_0(Y,k) \sim A^L(k)\cosh  
\left(\left(\alpha_M-\alpha(0)\right)^{\frac{1}{2}}Y\right)
\eeq
as $k\to \infty$ with $Y=O(1)^+$, and $A^L(k)$ a normalising factor to be determined. Similarly, for the boundary layer at $y=1$, we introduce,
\beq\elab{eqn4.4.104}
\bar{Y} = (y-1)k^{\frac{1}{2}} = O(1)^-
\eeq
as $k\to \infty$ in the boundary layer, after which we obtain,
\beq\elab{eqn4.4.105}
\psi_0(\bar{Y},k) \sim A^R(k)
\cosh
\left(
\left(\alpha_M-\alpha(0)\right)^{\frac{1}{2}}\bar{Y}\right)
\eeq
as $k\to \infty$ with $\bar{Y}=O(1)^-$, and $A^R(k)$ a normalising factor to be determined. The remaining asymptotic structure now has three regions, as follows:
\begin{itemize}
	\item region $L$ : $y\in \left(O(k^{-\frac{1}{2}})^+, y_M-O(k^{-\frac{1}{4}})^+\right)$
\item
\noindent region $R$: $y\in \left( y_M+O(k^{-\frac{1}{4}})^+, 1-O(k^{-\frac{1}{2}})^+\right)$
\item 
 spike region : $y\in \left(y_M-O(k^{-\frac{1}{4}})^+, y_M+O(k^{-\frac{1}{4}})^+\right)$.
\end{itemize}

We first move to region $L$. {The form of solution  \eref{eqn4.103} in the boundary layer  at $y=0$}, leads to a WKB form for the solution in this region, which, after asymptotic matching with \eref{eqn4.103}, gives,
\beq\elab{eqn4.106}
\psi_0(y,k) \sim A^L(k)\text{cosh}\left(k^{\frac{1}{2}}\int_{0}^{y}{(\alpha_M-\alpha(s))^{\frac{1}{2}}}ds\right)
\eeq
as $k\to \infty$ with $y\in \left(O(k^{-\frac{1}{2}})^+, y_M-O(k^{-\frac{1}{4}})^+\right)$. Similarly, in region $R$ we have,
\beq\elab{eqn4.107}
\psi_0(y,k) \sim A^R(k)\text{cosh}\left(k^{\frac{1}{2}}\int_{y}^{1}{(\alpha_M-\alpha(s))^{\frac{1}{2}}}ds\right)
\eeq
as $k\to \infty$ with $y\in \left( y_M+O(k^{-\frac{1}{4}})^+, 1-O(k^{-\frac{1}{2}})^+\right)$. We finally move to the spike region. In this region we introduce
\beq\elab{eqn4.108}
w = (y - y_M)k^{\frac{1}{4}} = O(1)
\eeq
as $k\to \infty$. The expansions \eref{eqn4.106} and \eref{eqn4.107} require us to write
\beq\elab{eqn4.109}
\psi_0(w,k) = A(k)\widetilde{\psi}_0(w) + o(A(k))
\eeq
as $k\to \infty$ with $w=O(1)$. Here $A(k)\to \infty$ as $k\to \infty$ is to be determined. Substitution from \eref{eqn4.108}, \eref{eqn4.109} and \eref{eqn4.100} into equation \eref{eqn4.73} gives, at leading order,
\begin{subequations}\elab{eqn4.1104.112}
\beq\elab{eqn4.110}
\widetilde{\psi}_{0ww} - \left(\frac{1}{2}(-\alpha_M'')w^2 - \widetilde{\lambda}\right) \widetilde{\psi}_0 = 0,~~w\in \mathbb{R},
\eeq
\beq\elab{eqn4.111}
\widetilde{\psi}_0(w)>0~~\forall~~w\in \mathbb{R},
\eeq
whilst matching with region $L$ and region $R$ (via Van Dyke's matching principle  \cite{VanDyke1975}) requires, on using \eref{eqn4.106} and \eref{eqn4.107},
\beq\elab{eqn4.112}
\widetilde{\psi}_0(w) \sim \text{exp}\left(-\frac{(-\alpha_M'')^{\frac{1}{2}}}{2^{\frac{3}{2}}}w^2\right)~~\text{as}~~|w|\to \infty,
\eeq
\end{subequations}
{and also that} $A^L(k)$, $A^R(k)$ and $A(k)$   satisfy the two conditions,
\beq\elab{eqn4.113}
 A^L(k)\text{exp}\left(k^{\frac{1}{2}}\int_{0}^{y_M}{(\alpha_M-\alpha(s))^{\frac{1}{2}}}ds\right)=A(k),
\eeq
\beq\elab{eqn4.114}
 A^R(k)\text{exp}\left(k^{\frac{1}{2}}\int_{y_M}^{1}{(\alpha_M-\alpha(s))^{\frac{1}{2}}}ds\right)=A(k).
\eeq
{In addition} the $L^1([0,1])$ normalising condition \eref{eqn4.80} requires,
\beq\elab{eqn4.115}
A(k)\int_{-\infty}^{\infty}{\widetilde{\psi}_0(w)}dw = k^{\frac{1}{4}}.
\eeq
In the above 
\beq
\alpha_M'' = \alpha''(y_M) < 0
\eeq
 restricting attention to a nondegenerate  interior maximum for $\alpha(y)$ (the case when the interior {maximum} is degenerate can be similarly addressed). Now, \eref{eqn4.110}-\eref{eqn4.112}
 is a linear self-adjoint singular Sturm-Liouville eigenvalue problem on the whole real line. The positivity  condition \eref{eqn4.111} determines that 
\emph{$\widetilde{\lambda}$ must be 
the smallest (principal) eigenvalue}.
  The condition \eref{eqn4.115} then normalises this eigenfunction. In fact, we can determine the lowest eigenvalue and the associated eigenfunction by inspection, which gives,
  \begin{subequations}
\beq\elab{eqn4.116}
\widetilde{\psi}_0(w) = \text{exp}\left(-\frac{(-\alpha_M'')^{\frac{1}{2}}}{2^{\frac{3}{2}}}w^2\right)~~\forall~~w\in \mathbb{R},
\eeq
with,
\beq\elab{eqn4.117}
\widetilde{\lambda} =\frac{1}{\sqrt{2}}(-\alpha_M'')^\frac{1}{2}.
\eeq
\end{subequations}
It then follows, via \eref{eqn4.113}--\eref{eqn4.116}, that,
\begin{subequations}
\beq\elab{eqn4.118}
A(k) = \frac{(-\alpha_M'')^{\frac{1}{4}}}{2^{\frac{3}{4}}\sqrt{\pi}}k^{\frac{1}{4}},
\eeq
\beq\elab{eqn4.119}
A^L(k) =  \frac{(-\alpha_M'')^{\frac{1}{4}}}{2^{\frac{3}{4}}\sqrt{\pi}}k^{\frac{1}{4}}\text{exp}\left(-k^{\frac{1}{2}}\int_{0}^{y_M}{(\alpha_M-\alpha(s))^{\frac{1}{2}}}ds\right),
\eeq
\beq\elab{eqn4.120}
A^R(k) =  \frac{(-\alpha_M'')^{\frac{1}{4}}}{2^{\frac{3}{4}}\sqrt{\pi}}k^{\frac{1}{4}}\text{exp}\left(-k^{\frac{1}{2}}\int_{y_M}^{1}{(\alpha_M-\alpha(s))^{\frac{1}{2}}}ds\right).
\eeq
\end{subequations}
Finally, via \eref{eqn4.100} and \eref{eqn4.117}, we have,
\beq\elab{eqn4.121}
\lambda_0(k) = -\alpha_Mk + \frac{1}{\sqrt{2}}(-\alpha_M'')^\frac{1}{2}k^{\frac{1}{2}} + O(1)
\eeq
as $k \to \infty$, and so \eref{eqn4.82} gives,
\beq\elab{eqn4.122}
v_1 = \alpha_M - \frac{1}{\sqrt{2}}(-\alpha_M'')^{\frac{1}{2}}k^{-\frac{1}{2}} + O(k^{-1})
\eeq
as $k\to \infty$. The analysis in this subcase is now complete. However, it is instructive to note that the structure will be adjusted when $y_M\to 0$ or $1$ as $k\to \infty$, and in particular, when $y_M=O(k^{-\frac{1}{4}})^+$ or $y_M=1 - O(k^{-\frac{1}{4}})^+$ as $k\to \infty$. A detailed consideration of the case when $y_M=O(k^{-\frac{1}{4}})^+$ as $k\to \infty$, reveals that the leading order form of $\lambda_0(k)$ remains unchanged, but the correction is influenced, so that, in these cases,
\beq
\lambda_0(k) = -\alpha_Mk + o(k)
\eeq
as $k\to \infty$. Evidently, the case of Couette flow falls into this category, and the details for this case are developed in \cite{HaynesVanneste2014a}, where it is established that the above correction is, in fact, of $O(k^{\frac{2}{3}})$. For brevity we do not pursue these special cases further.

\subsubsection{(SL) when $k=O(1)^+$}
In this case, details of (SL) must be obtained numerically. However, since we are dealing with the principal eigenvalue, whose normalised eigenfunction is strictly positive, a simple integration of equation \eref{eqn4.73} over the interval $[0,1]$,  using conditions \eref{eqn4.74}, and rearranging, establishes the following bounds on $\lambda_0(k)$, namely
\beq\elab{bounds}
-k\alpha_M < \lambda_0(k) < -k\alpha_m
\eeq
for all $k>0$, with
\beq\elab{bounds'}
\alpha_m = \inf_{y\in [0,1]}\alpha(y) < 0,
\eeq
when $\alpha(y)$ is nontrivial. We can thus conclude, from \eref{eqn4.82}, that,
\beq\elab{bounds'''}
\alpha_m < v_1(k) < \alpha_M
\eeq
for all $k>0$. In fact, we can improve the upper bound on $\lambda_0(k)$, by again using the strict positivity of the principal eigenfunction to establish,
after dividing through 
 equation \eref{eqn4.73} by $\psi$ and, using integration by parts together with \eref{eqn1.4}, that
\beq\elab{ubound}
\lambda_0(k) =   -\int_0^1{\frac{(\psi_0'(y,k))^2}{(\psi_0(y,k))^2}}dy < 0
\eeq
 for each $k>0$, when $\alpha(y)$ is nontrivial, and so we then  have the improved lower bound
\beq\elab{ubound'}
v_1(k) >0
\eeq
for all $k>0$.
 \begin{figure} [t] 
 		     \centering
 		    \includegraphics[width=0.5\textwidth]{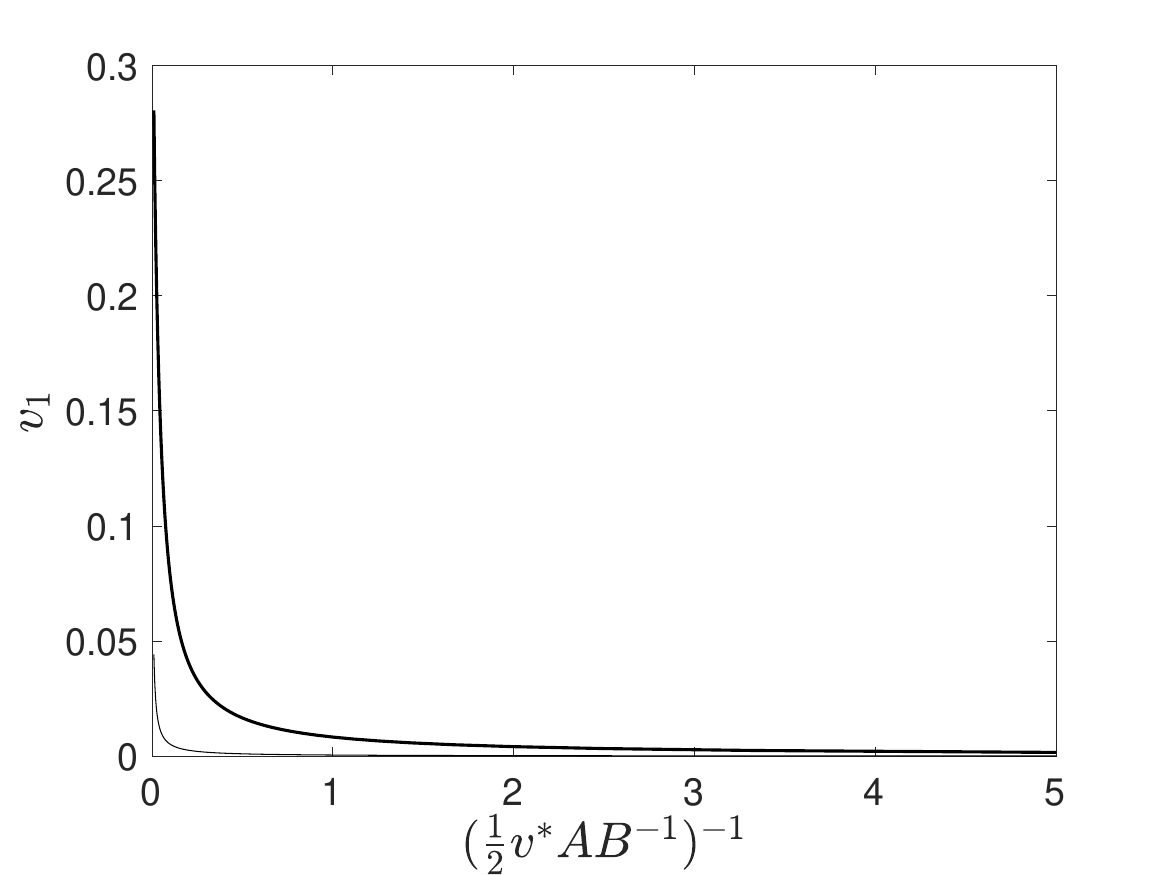}
 \caption{A graph of the scaled leading-order correction to the speed  as $A\to 0$ with $B=O(A)$ obtained  from (SL)  as a function of $k^{-1}=(\frac{1}{2}v^{\ast}AB^{-1})^{-1}$ for the  Couette flow (thick solid lines) and the Poiseuille flow (thin solid lines).  
} 
 \flab{speed_correction_Asmall}
 \end{figure} 
\medskip
It is instructive at this stage to summarise the asymptotic results concerning PTW solutions as $A\to 0$, the case of weak advection. We have considered the situation when $0<B\ll1$, $B=O(1)$ and $B\gg1$, in the natural distinguished limits, as $A\to 0$. We have shown that, when $A\ll1$, then for each $B>0$, there exists a unique PTW solution to (BVP) (correspondingly (QBVP)), with propagation speed denoted by $v=\hat{v}(A,B,u_c)$. With regards to this propagation speed, we have determined that, \emph{as $A\to 0$},
\begin{itemize}
\item[(i)] $B\gg1$ or $B=O(1)$
\beq\elab{eqn4.123}
\hat{v}(A,B,u_c) = v^{\ast}(u_c) + O\left(\frac{A^2}{B}\right)
\eeq
\item[(ii)] $B\ll1$ ($B=O(A)$)
\beq\elab{eqn4.125}
\hat{v}(A,B,u_c) = v^{\ast}(u_c) - Ak^{-1}\left(\frac{A}{B}\right)\lambda_0\left(k\left(\frac{A}{B}\right)\right) + O(A^2).
\eeq
Here
\beq\elab{eqn4.126}
k(\mu) = \frac{1}{2}v^{\ast}(u_c)\mu~~\forall~~\mu\in (0,\infty)
\eeq
  and $\lambda_0(k)$ is the principal eigenvalue of the Sturm-Liouville problem (SL). 
  We have established, upon use of results in \cite{HaynesVanneste2014a}, that when $k\to 0$ (so that   $0<A\ll B\ll 1$), then,
\beq\elab{eqn4.127}
\lambda_0(k) =-k^2\left\langle\left(\int_{0}^y\alpha(s)ds\right)^2\right\rangle+O(k^3)
\eeq
whilst when $k\to \infty$ (so that $0<B\ll A\ll 1$) then,
\beq\elab{eqn4.129}
k^{-1}\lambda_0(k) = -\alpha_M + \frac{1}{\sqrt{2}}(-\alpha_M'')^\frac{1}{2}k^{-\frac{1}{2}} + O(k^{-1}).
\eeq
\end{itemize}
We observe that the effects of weak advection on the PTW propagation speed is always of $o(1)$ as $A\to 0$, but becomes more significant with decreasing $B$. It is also worth noting here that expansion \eref{eqn4.125} forms a composite expansion for $\hat{v}(A,B,u_c)$ as $A\to 0$ for all $B>0$, with error uniformly of $O(A^2)$. In addition, it follows from this observation and the limiting forms \eref{eqn4.127} and \eref{eqn4.129}, that  the correction to the propagation speed is of $O(A)$ when $0<B \le O(A)$, initially decreasing from the value $A\alpha_M$ as $A^{1/2}B^{1/2}$, and continuing to decrease with increasing $B$, becoming of $O(A^2)$ for $B\ge O(1)$, and decreasing at a rate of $A^2B^{-1}$. We finally observe that the two complimentary asymptotic forms for the propagation speed as $A\to 0$, given in \eref{eqn4.123} and \eref{eqn4.125} match with each other according to Van Dyke's asymptotic matching principle \cite{VanDyke1975}.   
 To conclude this section, we now compare the moving boundary location in the PTW solutions. Assembling the results of this section we obtain, \emph{as $A\to 0$}, that the moving boundary is at spatial location $z=\zeta(y,A,B)$ for $y\in [0,1]$, with,
\begin{itemize}
    \item [(i)] $B\gg1$
\beq
\zeta(y,A,B) = O\left(\frac{A}{B}\right),~~y\in [0,1]
\eeq
\item [(ii)] $B=O(1)$
 \begin{subequations}
\beq
\zeta(y,A,B) = \frac{1}{B}\left(\langle\phi(\cdot)\rangle - 
\phi(y)\right)A + O(A^2),~~y\in [0,1]
\eeq
with
\beq
\phi(y)=\int_0^y\left({\int_0^s{\alpha(p)}dp}\right)ds.
\eeq
\end{subequations}
\item [(iii)]
 $B\ll1$ ($B=O(A)$)
\beq
\zeta(y,A,B) = \frac{2}{v^{\ast}(u_c)}(\log\psi_0(y,k) - \langle\log\psi_0(\cdot,k)\rangle) + O(A),~~y\in [0,1],
\eeq
with $k=\frac{1}{2}v^{\ast}(u_c)AB^{-1}$. Here $\psi_0:[0,1]\to \mathbb{R}$ is the principal, $L^1$-normalised eigenfunction of (SL). The structure of $\psi_0$ as $k\to 0$ ($A\ll B\ll 1$) and as $k\to \infty$ ($B\ll A\ll 1$) is as given in subsections 4.3.1 and 4.3.2. 
\end{itemize}

\begin{figure} [t] 
 		     \centering
 		   \begin{minipage}{0.49\textwidth}
 		    \includegraphics[width=\textwidth]{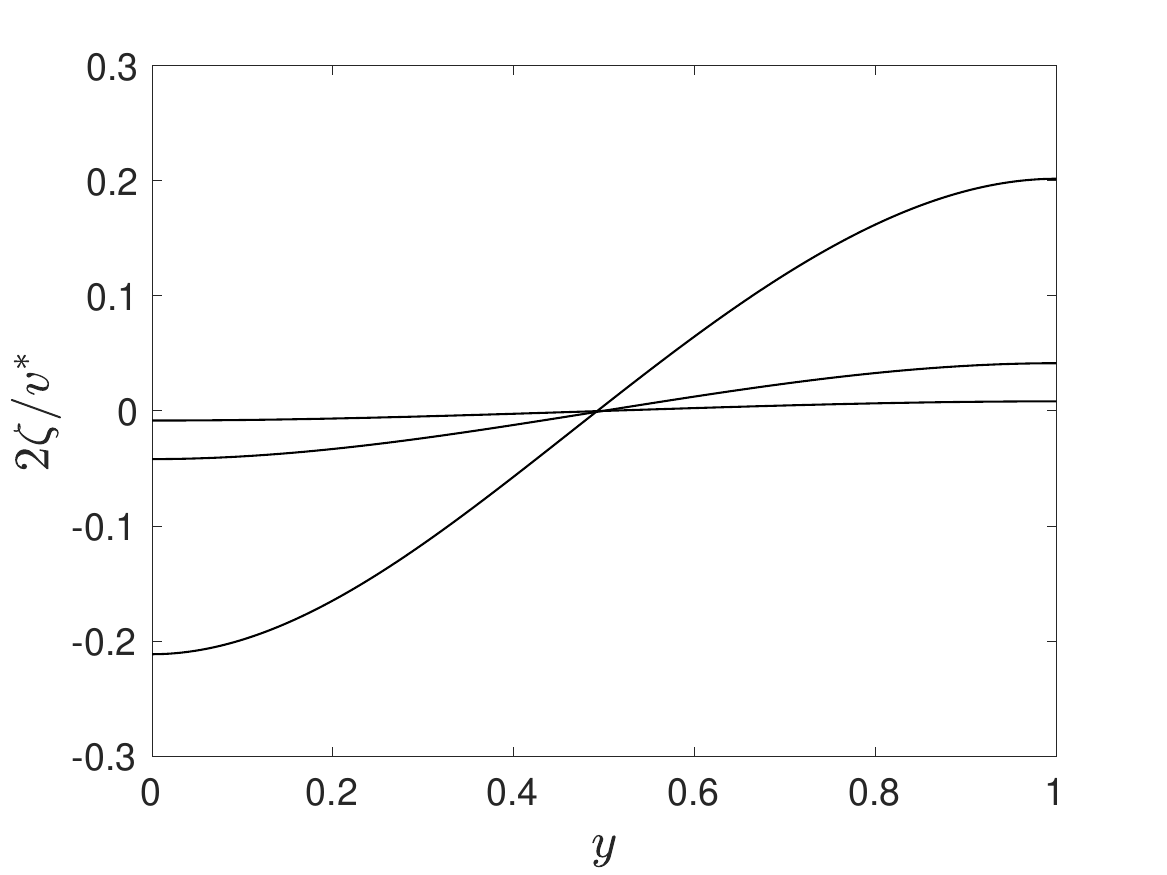}\\
 		     \end{minipage}
 		     \hfill
 		     \begin{minipage}{0.49\textwidth}
 		  	   \includegraphics[width=\textwidth]{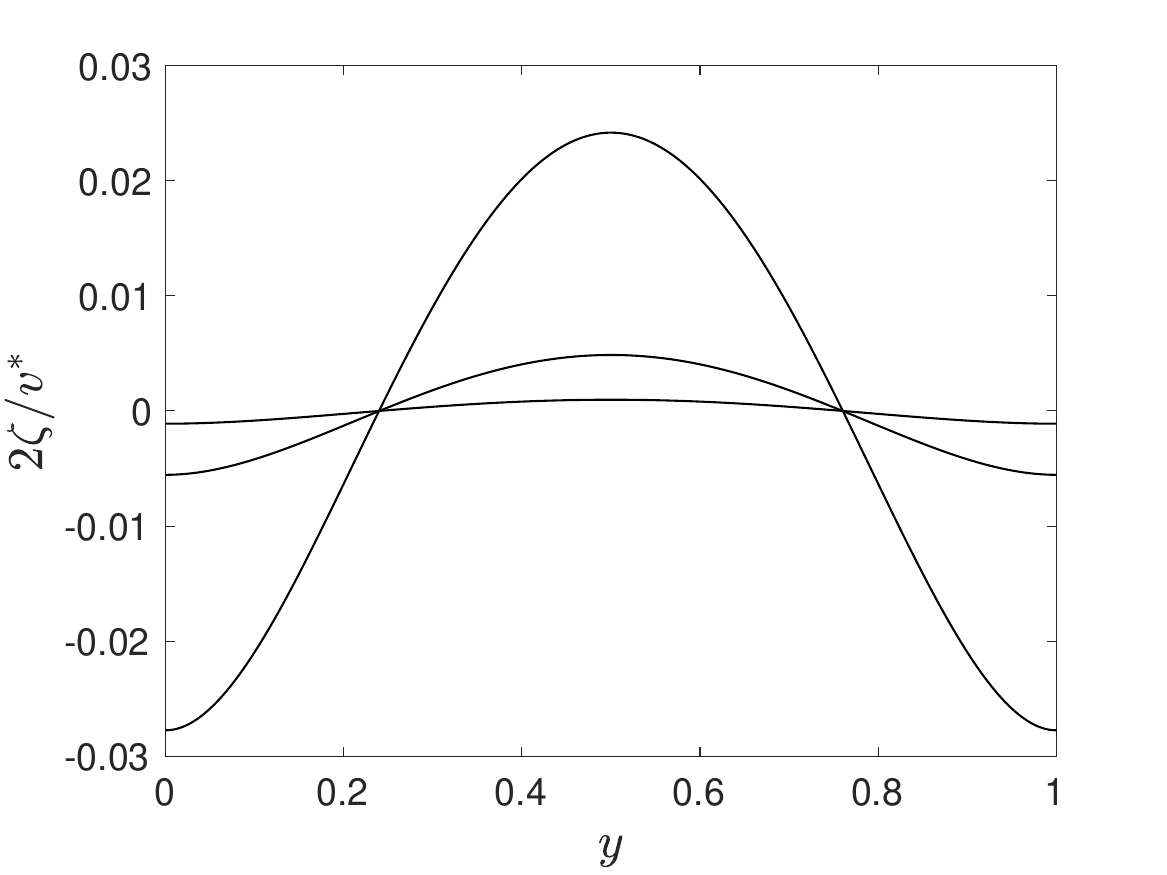}\\ 
			  \end{minipage}
 \caption{A graph of the scaled interface $2\zeta/v^*$ as $A\to 0$ with $B=O(A)$ for the  Couette flow (left) and the Poiseuille flow (right)  obtained by solving (SL)   for $k=1/5$, $k=1$ and $k=5$.  
The maximum of the interface decreases with {increasing} $k$.} 
 \flab{interface_correction_Asmall}
 \end{figure} 

Figure \fref{speed_correction_Asmall} shows the correction to speed of propagation of the PTW solution (divided by $A$) computed by numerical solution of the  (SL) eigenvalue problem 
for the Couette   and Poiseuille \eref{CouettePoiseuille} flows. Clearly, 
the enhancement to the speed is, for fixed (small) $A$, greater for the Couette flow than 
for the Poiseuille flow,  decreasing monotonically as either the front thickness or  cut-off value increase. 
Figure \fref{interface_correction_Asmall} shows the leading order scaled interface $2\zeta(y,A,B)/v^*(u_c)$, deduced from \eref{eqn4.81} for three different values of $k=\frac{1}{2}v^{\ast}(u_c)AB^{-1}$. The interface becomes increasingly flat as either  the front thickness or  cut-off value increase.

We now move on to consider PTW solutions when advection and streamwise diffusion are balanced.

\section{Slowly varying or rapidly varying front with balanced advection: $B\to \infty$ or $B\to 0$ with $A=O(1)$}
In this section we consider (IBVP) ((QIVP)) and (BVP) ((QBVP)) in the final case when the advective fluid velocity is comparable to the advectionless front propagation speed, so that $A=O(1)$. We consider this situation first when the channel width is small compared to the advectionless  front thickness ($B\gg1$), and then when the channel width is large compared to the advectionless front thickness ($B\ll1$). We begin with the former case.

\subsection{Slowly varying front: $B\to\infty$}
Following details very similar to those presented in section 3 (and without repetition) we have, for (IBVP),
\beq\elab{eqn5.1}
u(x,y,t;A,B) = \bar{u}(x,t) + O(B^{-1})
\eeq
as $B\to \infty$, with $(x,y,t)\in \mathbb{R}\times[0,1]\times[0,\infty)$, where $\bar{u}$ satisfies the one-dimensional evolution problem comprising the PDE \eref{eqn3.16}, \emph{with $\bar{B}^{-1}=0$}, subject to the initial and boundary conditions \eref{eqn3.8} and \eref{eqn3.9}. It then follows, regarding PTW solutions, that the propagation speed has
\beq\elab{eqn5.2}
\hat{v}(A,B,u_c) = v^{\ast}(u_c) + O(B^{-1})
\eeq
as $B\to \infty$ with $A=O(1)$. With the PTW moving boundary at $z=\zeta(y,A,B)$, it also follows from \eref{eqn5.1} that,
\beq\elab{eqn5.3}
\zeta(y,A,B) = O(B^{-1}),~~y\in [0,1]
\eeq
as $B\to \infty$ with $A=O(1)$. 

\subsection{Rapidly varying front: $B\to 0$}
We restrict attention to PTW solutions, and the most convenient formulation to work with is that given in (BVP). It is convenient to write
\beq\elab{eqn5.4}
 B=\widetilde{\eps}^2.
\eeq
The form of the PDE \eref{eqn2.53} with $A=O(1)$ as $\widetilde{\eps}\to 0$, leads us to write,
\beq\elab{eqn5.5}
U_T(z,y;\widetilde{\eps}) = u_0\left(\left(z-\frac{z_0(y)}{\widetilde{\eps}}\right)\psi(y)\right) + O(\widetilde{\eps})
\eeq
as $\widetilde{\eps}\to 0$ with $(z,y)\in \mathbb{R}\times [0,1]$, and $z_0,\psi:[0,1]\to \mathbb{R}$, $u_0:\mathbb{R}\to \mathbb{R}$ each having $z_0$, $\psi$, $u_0$ $=O(1)$ as $\widetilde{\eps}\to 0$. In addition we require,
\beq\elab{eqn5.6}
\psi(y)>0~~\forall~~y\in [0,1],
\eeq
via \eref{eqn2.57}. In addition to \eref{eqn5.5} we also expand,
\beq\elab{eqn5.7}
v(\widetilde{\eps}) = v_0 + O(\widetilde{\eps})
\eeq
as $\widetilde{\eps})$. On substitution from \eref{eqn5.5} and \eref{eqn5.7} into (BVP) we obtain, at leading order,
\beq\elab{eqn5.8}
\begin{split}
 (
 \psi^2(y)(1 + z_0'(y)^2) + 
 & 2z_0(y)z_0'(y)\psi(y)\psi'(y) + z_0^2(y)\psi'(y)^2
 )u_0''(\xi)\\ 
 &+ \psi(y)(v_0-A\alpha(y))u_0'(\xi) + f_c(u_0(\xi)) = 0,~~(\xi,y)\in \mathbb{R}\times(0,1),
\end{split}
\eeq
with $'$ representing differentiation in the respective arguement, and where $\xi$, replacing $z$, is given by,
\beq\elab{eqn5.8'}
\xi = \left(z - \frac{z_0(y)}{\widetilde{\eps}}\right)\psi(y).
\eeq
Since $\xi$ and $y$ are independent coordinates, the solubility of \eref{eqn5.8} requires, without loss of generality,
\begin{subequations}\elab{eqn5.8''eqn5.10}
\beq\elab{eqn5.8''}
\psi^2(y)(1 + z_0'(y)^2) + 2z_0(y)z_0'(y)\psi(y)\psi'(y) + z_0^2(y)\psi'(y)^2\\ = 1,
\eeq
\beq\elab{eqn5.9}
\psi(y)(v_0 - A\alpha(y)) = v^{\ast}(u_c),
\eeq
for $y\in (0,1)$, on normalising $\psi$ so that,
\beq\elab{norm}
\sup_{y\in[0,1]}\psi(y)=1.
\eeq
The boundary conditions in (BVP) require,
\beq\elab{eqn5.10}
\psi'(0)=\psi'(1)=z_0'(0)=z_0'(1)=0.
\eeq
\end{subequations}
Before considering \eref{eqn5.8''eqn5.10} in more detail, we first write down the problem for $u_0(\xi)$, which becomes, via \eref{eqn5.8''eqn5.10}, and \eref{eqn5.6},
\begin{subequations}
\beq
u_0''(\xi) + v^{\ast}(u_c)u_0'(\xi) + f_c(u_0(\xi)) = 0,~~\xi\in \mathbb{R},
\eeq
subject to
\beq
u_0(\xi)\ge0~~\forall~~\xi\in \mathbb{R},
\eeq
\beq
u_0(\xi)\to
\begin{cases}
0~~\text{as}~~\xi\to \infty, \\
1~~\text{as}~~\xi\to {-\infty}.
\end{cases}
\eeq
\end{subequations}
This problem has the solution (as noted in previous sections, and see also \cite{Tisbury_etal2020a}),
\beq\elab{eqn5.11}
u_0(\xi(z,y)) = U_T(\xi(z,y)),~~(z,y)\in \mathbb{R}\times[0,1],
\eeq
and we recall from \cite{Tisbury_etal2020a} that,
\beq\elab{eqn5.12}
U_T'(s)<0~~\forall~~s\in \mathbb{R}.
\eeq
Now, following \eref{eqn5.5} and \eref{eqn5.11}, the moving boundary, where $U_T=u_c$, occurs at $\xi(z,y)=0$, and so, via \eref{eqn5.7}, when
\beq\elab{eqn5.13}
z = \frac{z_0(y)}{\widetilde{\eps}},~~y\in [0,1]
\eeq
and so, we must have, via (BVP) (\eref{eqn2.64}),
\beq\elab{eqn5.14}
\int_0^1{z_0(y)}dy = 0.
\eeq
To proceed we will limit attention to shear flows $\alpha(y)$ which have a single point of maximum velocity $\alpha_M$, occuring at $y=y_M$, where, without loss of generality, we may take $y_M\in [\frac{1}{2},1]$. Next we observe from \eref{eqn5.9} that the continuous functions $\psi(y)$ and $\alpha(y)$ achieve their maximum values on $[0,1]$, of unity and $\alpha_M$ respectively, at the same point $y=y_M$,
after which \eref{eqn5.9} evaluated at $y=y_M$, requires,
\beq\elab{eqn5.21}
v_0 = v^{\ast}(u_c) + A\alpha_M.
\eeq
Equation \eref{eqn5.9}   then becomes
\beq\elab{eqn5.22}
\psi(y) = v^{\ast}(u_c)(v^{\ast}(u_c) + A(\alpha_M - \alpha(y)))^{-1}~~\forall~~y\in [0,1].
\eeq
 \begin{figure} [t] 
 		     \centering
 		     \begin{minipage}{0.49\textwidth}
 		    \includegraphics[width=\textwidth]{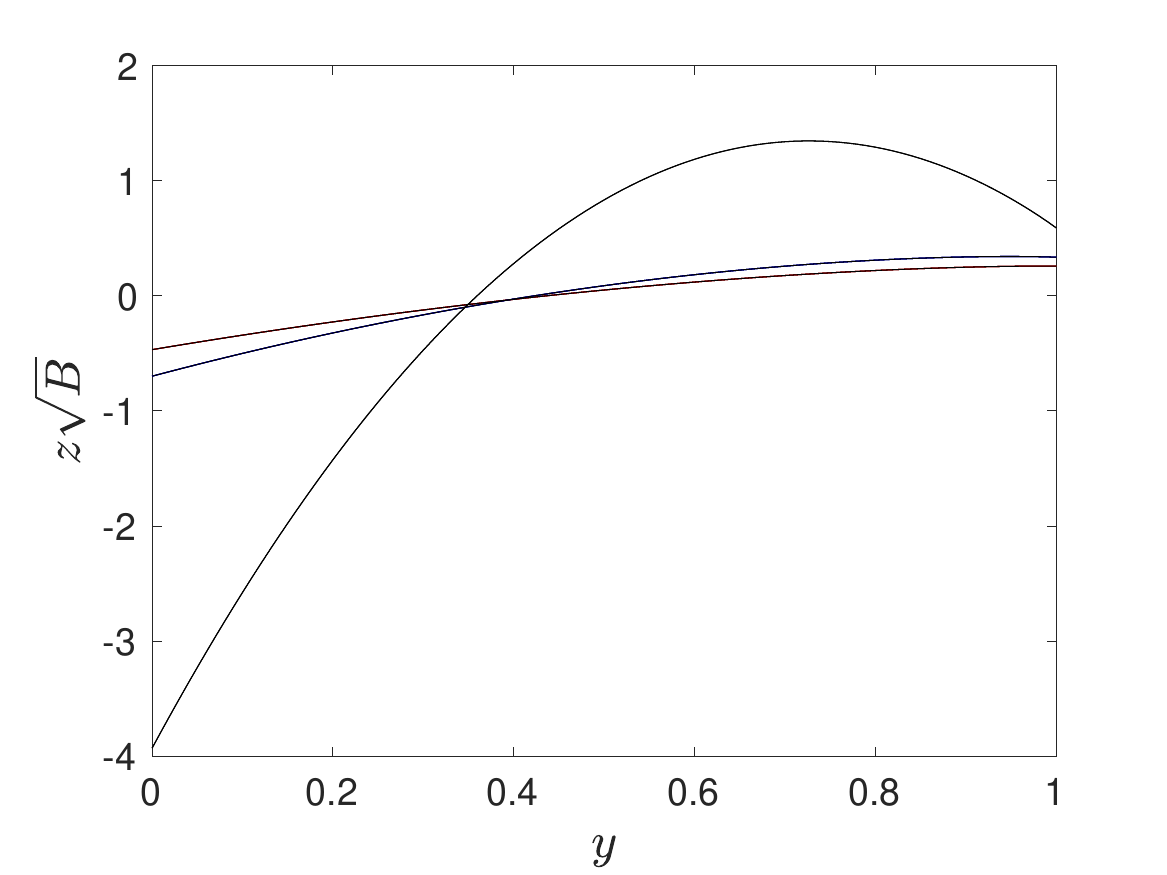}\\
 		     \end{minipage}
 		     \hfill
 		     \begin{minipage}{0.49\textwidth}
 		  	   \includegraphics[width=\textwidth]{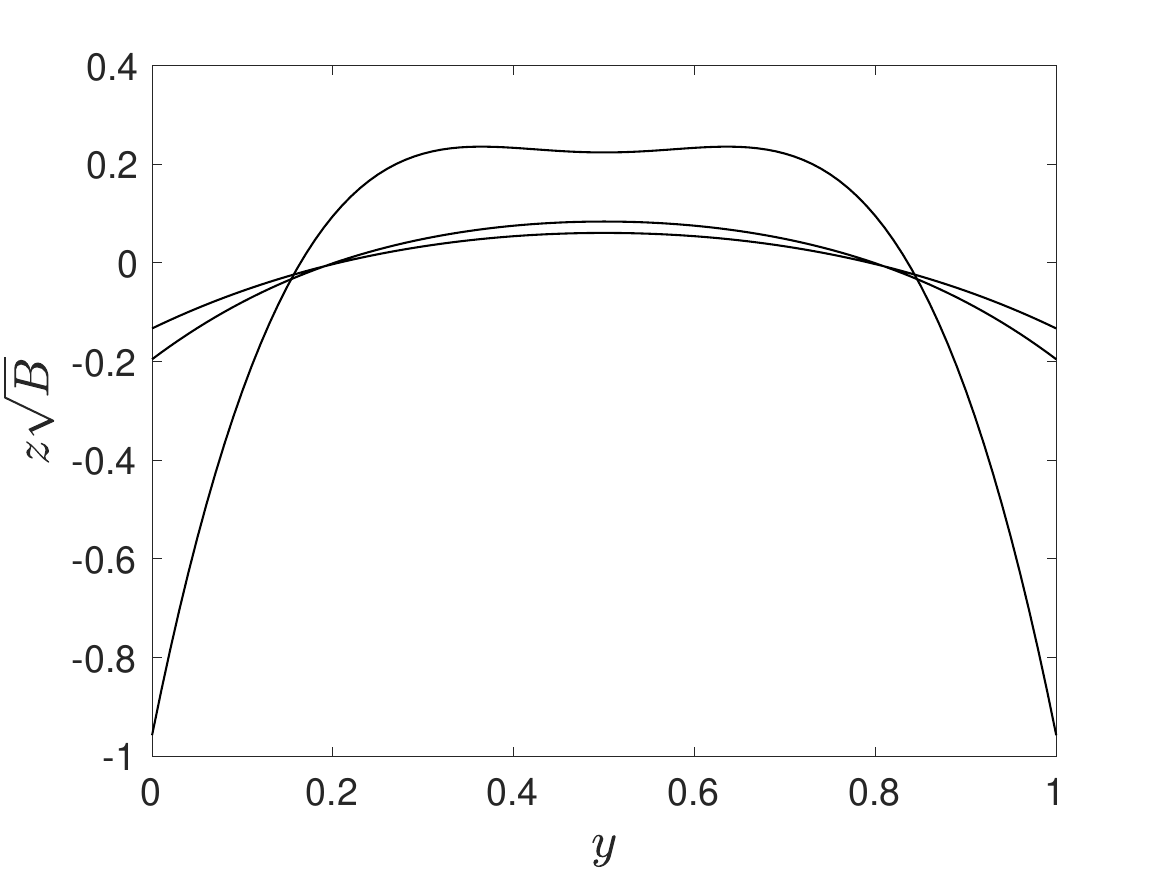}\\ 
			  \end{minipage}
 \caption{A graph of the scaled interface $z\sqrt{B}$ as $B\to 0$  for $A=1$ for the (left) Couette flow \eref{Couette} and (right) Poiseuille flow \eref{Poiseuille}. 
 These are numerically obtained  
 using \eref{eqn5.23} and \eref{eqn5.14} for $v^*(u_c)=0.1$, $1$ and $1.9$. 
The   maximum of the interface decreases with {increasing} $v^*(u_c)$.}

 \flab{z0_Poiseuille}
 \end{figure} 
We are left with \eref{eqn5.8} to determine $z_0(y)$. We also note that, in general, \eref{eqn5.22} will not satisfy the boundary conditions in \eref{eqn5.10}, and passive boundary layers will be needed, when $y=O(\widetilde{\eps}^{\frac{1}{2}})^+$ and $y=1-O(\widetilde{\eps}^{\frac{1}{2}})^+$. For brevity, we do not consider these boundary layers here. 
We observe from \eref{eqn5.22} that,
\beq\elab{eqn5.22'}
0<\psi(y) \le 1~~\forall~~y\in [0,1]
\eeq
as specified in the normalisation \eref{norm}. In addition the stationary points of $\alpha(y)$ and $\psi(y)$ are coincident. Now, equation \eref{eqn5.8''} can be re-written as
\beq
\left((\psi(y)z_0(y))'\right)^2 + \psi^2(y) - 1 = 0,~~y\in (0,1).
\eeq
Therefore,
\beq
(\psi(y)z_0(y))' = 
\begin{cases}
\left(1-\psi^2(y)\right)^{\frac{1}{2}},~~0\le y<y_M,\\
-\left(1-\psi^2(y)\right)^{\frac{1}{2}},~~y_M\le y\le 1,
\end{cases}
\eeq
so that an integration gives,
\beq\elab{eqn5.23}
z_0(y) = c\psi^{-1}(y) + \psi^{-1}(y)
\begin{cases}
\int_0^y{\left(1-\psi^2(s)\right)^{\frac{1}{2}}}ds,~~0\le y<y_M,\\
\int_0^{y_M}{\left(1-\psi^2(y)\right)^{\frac{1}{2}}}ds - \int_{y_M}^y{\left(1-\psi^2(y)\right)^{\frac{1}{2}}}ds,~~y_M\le y \le1,
\end{cases}
\eeq
with the constant $c$ determined by the condition \eref{eqn5.14}. We note that, as required, $z_0\in C^1([0,1])$, with passive boundary layers required when $y, (y-1) = O(\widetilde{\eps}^{\frac{1}{2}})$ over which the boundary conditions \eref{eqn5.10} will be satisfied. In summary we have, in this case of balanced advection, with the channel width large compared to the streamwise diffusion length scale (based on the reaction time scale), a PTW solution for each $A>0$, with propagation speed
\beq
\hat{v}(A,B,u_c) = v^{\ast}(u_c) + A\alpha_M + O(B^{\frac{1}{2}})
\eeq
as $B\to 0$ with $A=O(1)$. In addition, the moving boundary has location,
\beq
z = B^{-\frac{1}{2}}z_0(y) + O(1),~~y\in [0,1],
\eeq
as $B\to 0$ with $A=O(1)$, $z_0:[0,1]\to \mathbb{R}$ given by \eref{eqn5.23} and $\psi:[0,1]\to \mathbb{R}$ given by \eref{eqn5.9}.
We now compute the interface deformation for the Poiseuille flow \eref{Poiseuille} by numerically integrating  \eref{eqn5.23} with  \eref{eqn5.14} 
   for three different values of $v^*(u_c)$ corresponding to small, intermediate and large cut-off $u_c$.
As the cut-off increases, the maximum of the interface and its deformation increases.

\section{Conclusions}
In this paper we have considered  a  reaction-diffusion process evolving inside an infinite channel in the presence of a shear flow. 
The reaction function is of standard KPP-type, but experiences a cut-off in the reaction rate below the normalised cut-off concentration $u_c\in(0,1)$. 
We have  formulated this initial-value boundary problem (IBVP) in terms of  
{an equivalent moving boundary evolution problem} 
  (QIVP) and examined both of these problems in  a number of natural asymptotic limits relating to the two  non-dimensional (positive) parameters $A$ and $B$ describing the strength of the flow and front thickness, whilst keeping $u_c$ fixed. We have  established in all cases considered, that a unique (up to spatial translation) PTW solution $U_T(x-v t,y)=u(x,y,t)$ exists, and in a number of cases,    that this forms a large-$t$ attractor for (IBVP). 
We used the method of matched asymptotic expansions to  determine the  detailed structure of the PTW solutions and their speed of propagation and, in particular, their dependence on the parameters. 
It is of interest to examine the nature of 
the propagation speed $\hat{v}(A,B,u_c)$  
in the three complementary regimes $A\to\infty$,
$A=O(1)$ and $A\to 0$ over the full range of 
$B$.  
Using the results captured by Figure \fref{speed_correction_Asmall},  
  $\hat{v}(A,B,u_c)$ is fully determined, up to the stated orders in $A$, as a function of $B$ as $A\to 0$. Partial results obtained in the other two regimes are completed by extracting and extrapolating limiting forms of the available asymptotic expressions. 
 We therefore expect that for $A=O(1)$ and $A\to \infty$, the speed of propagation decreases monotonically with $B$.
  We anticipate that the approach developed in this paper will be readily adaptable to the case of 
  two corresponding problems, when the KPP-type cut-off reaction function is replaced by a broader class of cut-off reaction functions 
  or when the shear flow is unsteady, varying periodically with  time. 
Finally, it would be interesting to determine whether   similar qualitative effects arise for a certain stochastically perturbed KPP equation
obtained from an alternative
 model whose purpose is also to account for microscopic discrete particles \cite{ConlonDoering2005,DoeringTesser2014}. 
 Mueller, Mytnik and Quastel \cite{Mueller_etal2011} have shown that in the absence of a shear, 
 the difference between the speed obtained from this model and $v^*(u_c)$
  obtained from the deterministic cut-off model considered here is small when the   cut-off $u_c$ and noise are both small. Whether this difference continues to be small in the presence of a shear flow is presently unknown.

\medskip

\noindent
\textbf{Acknowledgments.}
The authors  acknowledge A. D. O. Tisbury for his early contribution on this project. We also thank the referees for their constructive comments, which have led to an overall improvement in the presentation of the paper.

\appendix

\section{{Asymptotic solution to (QIVP) as $t\to 0^+$}}

 {As $t\to 0^+,$ 
 the nature of the discontinuity in the initial condition of (QIVP) and the requirement for this to be locally spatially smoothed in the small time limit, demands a diffusion balance at leading order in (QIVP). As a consequence of this,} it follows that   
  the solution to (QIVP) develops in four principal asymptotic regions, which are,
\begin{itemize}

\item Region $\textrm{I}_L$ : $\xi'=O(t^{\frac{1}{2}})^-$,  $y'=O(1)$ with $u=O(1)$ as $t\to 0^+$

\item Region $\textrm{I}_R$ : $\xi'=O(t^{\frac{1}{2}})^+$,   $y'=O(1)$ with $u=O(1)$ as $t\to 0^+$

\item Region $\textrm{II}_L$ : $\xi'=O(1)^-$,   $y'=O(1)$ with $u=1-o(1)$ as $t\to 0^+$

\item Region $\textrm{II}_R$ : $\xi'=O(1)^+$,   $y'=O(1)$ with $u=o(1)$ as $t\to 0^+$
\end{itemize}
with, in addition,
\beq\elab{eqn2.37}
s(t)=o(1)~~,~~\zeta(y',t)=o(1)~~\text{as}~~t\to 0^+.
\eeq
{A similar structure characterises the case in the absence of advection and  this is developed in \cite{Tisbury_etal2020b}, where a qualitative sketch of the above regions is also included.}
Throughout $O(\cdot)^{\pm}$ denotes $O(\cdot)$ being positive or negative, respectively. A balance of terms in the PDE \eref{eqn2.10} indicates that, with regions $\textrm{I}_{L,R}$ having $\xi=O(t^{\frac{1}{2}})$ as $t\to 0^+$, then,
\[
s(t)=O(t^{\frac{1}{2}}),~~\zeta(y',t)=O(t^{\frac{1}{2}})
\]
as $t\to 0^+$. Thus we expand in the form,
\beq\elab{eqn2.392.40}
s(t) = s_0t^{\frac{1}{2}} + s_1t + O(t^{\frac{3}{2}})
\quad
\text{and}
\quad 
\zeta(y',t) = \zeta_0(y')t^{\frac{1}{2}} + \zeta_1(y')t + O(t^{\frac{3}{2}})
\eeq
as $t\to 0^+$, with $y'\in [0,L]$. 
The details in regions $\textrm{I}_{L,R}$ are straightforward, obtained by introducing the coordinate $\eta=\xi 't^{-1/2}$ 
and in both regions expanding 
\beq\elab{smalltexpansion}
u(\eta,y',t)=u_{0}(\eta,y')+t^{1/2}u_{1}(\eta,y')+O(t),\quad\text{as}\quad t\to 0^+.
\eeq
Substituting expansions \eref{eqn2.392.40},   and \eref{smalltexpansion}
into the equation for $u$ we obtain that  $u_0=u_0(\eta)$ and satisfies
\begin{subequations}
\beq\elab{leadingordersmallt}
u_{0\eta\eta}+\frac{1}{2}(\eta+s_0)u_{0\eta}=0,~~ \eta\in \mathbb{R} \setminus \{0\}
\eeq
  subject to 
\beq
u_0(0)=u_c, \quad
u_{0\eta}(0^-)=u_{0\eta}(0^+)
\eeq
together with matching conditions with regions II$_L$  and  II$_R$   which requires,
\beq
u_0(\eta)\to
\begin{cases}
1\quad\text{as}\quad \eta\to-\infty,\\
0\quad\text{as}\quad \eta\to\infty.
\end{cases}
\eeq
\end{subequations}
It is straightforward to determine
\beq\elab{u0smallt}
u_0(\eta)=\frac{1}{2}\text{erfc}(\frac{1}{2}(\eta+s_0)),
\quad
\text{where}
\quad 
s_0=2\text{erfc}^{-1}(2u_c)
\quad
\text{with}
\quad
\zeta_0(y')=0,
\eeq
for $(\eta,y')\in \mathbb{R}\times[0,L]$. 
At the next order, we obtain that 
$u_1=u_1(\eta,y')$ satisfies
\begin{subequations}\elab{nextordersmallt}
\beq
u_{1\eta\eta}+\frac{1}{2}(\eta+s_0)u_{1\eta}=A\bar\alpha(y'),~~(\eta,y')\in (\mathbb{R} \setminus \{0\})\times(0,L)
\eeq
subject to
\beq
u_1(0,y')=-\zeta_1(y')u_{0\eta}(0), \quad
u_{1\eta}(0^-,y')=u_{1\eta}(0^+,y'),
\eeq
with $y'\in [0,L]$,
together with matching conditions with regions II$_L$  and  II$_R$   which requires,
\beq
u_1(\eta,y)\to 0\quad\text{as}\quad|\eta|\to\infty ~~\text{uniformly for}~~y\in [0,L].
\eeq
\end{subequations}
Solving \eref{nextordersmallt}
yields 
\beq\elab{u1smallt}
u_1(\eta,y')=\frac{\zeta_1(y')}{2\sqrt{\pi}}e^{-\frac{1}{4}(\eta+s_0)^2}, 
\quad
\text{where}
\quad 
s_1=0, 
\quad
\text{with}
\quad
\zeta_1(y')=A\bar\alpha(y'),
\eeq
for $(\eta,y')\in \mathbb{R} \times[0,L]$. 
Thus, we have via \eref{u0smallt} and \eref{u1smallt} that
\begin{subequations}
\beq\elab{eqn2.41}
u(\eta,y',t) = \frac{1}{2}\text{erfc}(\frac{1}{2}(\eta+s_0)) + \frac{A\bar{\alpha}(y')}{2\sqrt{\pi}}e^{-\frac{1}{4}(\eta+s_0)^2}t^{\frac{1}{2}} + O(t)
\eeq
as $t\to 0^+$, uniformly for $(\eta,y')=(O(1)^{\pm})\times [0,L]$, with 
\beq\elab{eqn2.44}
s(t) = 2\text{erfc}^{-1}(2u_c)t^{\frac{1}{2}} + O(t^{\frac{3}{2}})
\eeq
and
\beq\elab{eqn2.45}
\zeta(y',t) = A\bar{\alpha}(y')t + O(t^{\frac{3}{2}})
\eeq
uniformly for $y'\in [0,L]$ as $t\to 0^+$.
\end{subequations}
An interesting interpretation of these asymptotic forms is that $s(t)$ gives the diffusive displacement, whilst $\zeta(y',t)$ gives the advective displacement, respectively, of the spatial moving boundary in the early times. However, reaction plays no role, up to $O(t^{\frac{3}{2}})$.

We now move on to regions $\mathrm{II}_{L,R}$. It follows from \eref{eqn2.41} that $u$ and $u-1$ are exponentially small in $t$ as $t\to 0^+$ in regions $\mathrm{II}_R$ and $\mathrm{II}_L$, respectively. The details are lengthy, but standard (see, for example,  \cite{Tisbury_etal2020b} for more details for the case when $\bar\alpha(y')=0$ for $y'\in[0,L]$), and are omitted for brevity. After matching with regions $\mathrm{I}_L$ and $\mathrm{I}_R$ respectively, we have  
\beq\elab{eqn2.46}
u(\xi',y',t) = H(-\xi') + \text{sgn}(\xi')\exp{(-\phi(\xi',y',t)t^{-1})}
\eeq
with $H(\cdot)$ being the usual Heaviside function, and
\beq\elab{eqn2.47}
\phi(\xi',y',t) = 
\frac{1}{4}\xi'^2 + \frac{1}{2}s_0\xi't^{\frac{1}{2}} - \left(\frac{1}{2}\log{t} - \frac{1}{4}s_0^2 - \frac{1}{2}\log{\pi} + \frac{1}{2}A\bar{\alpha}(y')\xi'- \log(|\xi'|)\right)t + O(t^{\frac{3}{2}})
\eeq 
as $t\to 0^+$, with $(\xi',y')\in (O(1)^{\pm})\times [0,L]$.

At this stage the asymptotic structure as $t\to 0^+$ is not quite complete. In obtaining the expansions in regions $\mathrm{I}_{L,R}$ and $\mathrm{II}_{L,R}$, we have been forced to neglect the Neumann boundary conditions at $y'=0$ and $y'=L$, and an examination of these expansions reveals that,
\beq\elab{eqn2.50}
u_{y'}\sim O(t^{\frac{1}{2}})   
\eeq
when $y'=0,L$ in regions $\mathrm{I}_{L,R}$ as $t\to 0^+$, whilst,
\beq\elab{eqn2.51}
u_{y'}\sim O(e^{-\phi_{L,R}t^{-1}})
\eeq
when $y'=0,L$ in regions $\mathrm{II}_{L,R}$ as $t\to 0^+$. Thus, \emph{weak} and \emph{passive} boundary layers are required when $y'=O(t^{\frac{1}{2}})$ and $y'=L-O(t^{\frac{1}{2}})$ as $t\to 0^+$, adjacent to both regions $\mathrm{I}_{L,R}$ and $\mathrm{II}_{L,R}$. These boundary layers are readily dealt with, and since they are passive in nature, for brevity we do not present details here. This completes the asymptotic structure of the solution to (QIVP) as $t\to 0^+$.

\bibliographystyle{abbrv}

\begin{thebibliography}{10}

 \bibitem{Berestycki2003}
 {\sc H.~Berestycki}, {\em The influence of advection on the propagation of
   fronts in reaction-diffusion equations}, in Nonlinear PDEs in Condensed
   Matter and Reactive Flows, H.~Berestycki and Y.~Pomeau, eds., vol.~569 of
   NATO Science Series C, Kluwer, Doordrecht, 2003.

 \bibitem{BerestyckiNirenberg1992}
 {\sc H.~Berestycki and L.~Nirenberg}, {\em Travelling fronts in cyclinders},
   Ann. I. H. Poincar\'e, 9 (1992), pp.~497--572.



 \bibitem{BrunetDerrida1997}
 {\sc E.~Brunet and B.~Derrida}, {\em Shift in the velocity of a front due to a
   cut-off}, Phys. Rev. E., 56 (1997), pp.~2597 -- 2604.

 \bibitem{Camassa_etal2010}
 {\sc R.~Camassa, Z.~Lin, and R.~M. McLaughlin}, {\em The exact evolution of the
   scalar variance in pipe and channel flow}, Commun. Math. Sci., 8 (2010),
   pp.~601--626.

 \bibitem{CoddingtonLevinson}
 {\sc E.~A. Coddington and N.~Levinson}, {\em Theory of ordinary differential
   equations}, McGraw-Hill, New York, 1955.

 \bibitem{ConlonDoering2005}
 {\sc J.~C. Conlon and C.~D. Doering}, {\em On travelling waves for the
   stochastic {F}isher-{K}olmogorov-{P}etrovskii-{P}iscounov equation}, 120
   (2005), pp.~421--477.

 \bibitem{DoeringTesser2014}
 {\sc C.~R. Doering and F.~Tesser}, {\em Discrete and Continuum Dynamics of
   Reacting and Interacting Individuals. In: Muntean A., Toschi F. (eds)
   Collective Dynamics from Bacteria to Crowds. CISM International Centre for
   Mechanical Sciences}, vol.~553, Springer, Vienna, 2014.

 \bibitem{Dumortier_etal2007}
 {\sc F.~Dumortier, N.~Popovic, and T.~J. Kaper}, {\em The critical wave speed
   for the {F}isher-{K}olmogorov-{P}etrovskii-{P}iscounov equation with
   cut-off}, Nonlinearity, 20 (2007), pp.~855--877.
   
 \bibitem{GilbargTrudinger2001}
 {\sc D.~Gilbarg, and  N.~S.~Trudinger}, {\em Elliptic Partial Differential Equations of Second Order}, Berlin-Heidelberg-New York-Tokyo, Springer-Verlag (1983).	



 \bibitem{Embid_etal1995}
 {\sc P.~F. Embid, A.~J. Majda, and P.~E. Souganidis}, {\em Comparison of
   turbulent flame speeds from complete averaging and the {G}‐equation}, Phys.
   Fluids, 7 (1997), pp.~2052--2060.

 \bibitem{Fisher1937}
 {\sc R.~A. Fisher}, {\em The wave of advance of advantageous genes}, Ann.
   Eugenics, 7 (1937), pp.~355--369.

 \bibitem{Freidlin1985}
 {\sc M.~I. Freidlin}, {\em Functional Integration and Partial Differential
   Equations}, Princeton University Press, 1985.

 \bibitem{GertnerFreidlin1979}
 {\sc J.~G{\"{a}}rtner and M.~I. Freidlin}, {\em On the propagation of
   concentration waves in periodic and random media.}, Soviet Math. Dokl., 20
   (1979), pp.~1282--1286.

 \bibitem{HaynesVanneste2014a}
 {\sc P.~H. Haynes and J.~Vanneste}, {\em Dispersion in the large-deviation
   regime. {P}art 1: shear flows and periodic flows}, J.\ Fluid Mech., 745
   (2014), pp.~321--350.

 \bibitem{Kolmogorov_etal1937}
 {\sc A.~N. Kolmogorov, I.~G. Petrovsky, and N.~S. Piskunov}, {\em {\'{E}}tude
   de l'{\'{e}}quation de la diffusion avec croissance de la quantit{\'{e}} de
   mati{\`{e}}re et son application {\`{a}} un probl{\`{e}}me biologique}, Bull.
   Univ. Moskov. Ser. Internat. Sect., 1 (1937), pp.~1--25.

 \bibitem{MajdaKramer1999}
 {\sc A.~J. Majda and P.~R. Kramer}, {\em Simplified models for turbulent
   diffusion: Theory, numerical modelling, and physical phenomena}, Phys.\ Rep.,
   314 (1999), pp.~237 -- 574.

 \bibitem{MajdaSouganidis1993}
 {\sc A.~J. Majda and O.~E. Sougadinis}, {\em Large-scale front dynamics for
   turbulent reaction--diffusion equations with separated velocity scales},
   Nonlinearity, 7 (1993), pp.~1--30.

 \bibitem{MallordyRoquejoffre1995}
 {\sc J.-F. Mallordy and J.-M. Roquejoffre}, {\em A parabolic equation of the
   {K}{P}{P} type in higher dimensions}, SIAM J. Math. Anal., 26 (1995),
   pp.~1--20.

\bibitem{Miller}
{\sc P.~D.~Miller}, {\em Applied Asymptotic Analysis}, Volume 75 of Graduate studies in Mathematics, American Mathematical Society, Providence RI, 2006. 

 \bibitem{Mueller_etal2011}
 {\sc C.~Mueller, L.~Mytnik, and J.~Quastel}, {\em Effect of noise on front
   propagation in reaction-diffusion equations of {KPP} type}, Invent. Math.,
   184 (2011), pp.~405--453.

 \bibitem{PapanicolaouXin1991}
 {\sc G.~Papanicolaou and J.~Xin}, {\em Reaction-diffusion fronts in
   periodically layered media}, J. Stat. Phys.,  (1991), pp.~915--931.

 \bibitem{Roquejoffre1997}
 {\sc J.~M. Roquejoffre}, {\em Eventual monotonicity and convergence to
   travelling fronts for the solutions of parabolic equations in cylinders},
   Ann. I. H. Poincar\'e, 14 (1997), pp.~499--552.

 \bibitem{Heinze_etal2001}
 {\sc A.~Stevens, G.~Papanicolaou, and S.~Heinze}, {\em Variational principles
   for propagation speeds in inhomogeneous media}, SIAM J. Appl. Math., 62
   (2001), pp.~129--148.

 \bibitem{Tisbury_etal2020a}
 {\sc A.~D.~O. Tisbury, D.~J. Needham, and A.~Tzella}, {\em The evolution of
   traveling waves in a {K}{P}{P} reaction-diffusion model with cut-off reaction
   rate. {I}. {P}ermanent form traveling waves}, Stud. in App. Math., 146
   (2021), pp.~301--329.

 \bibitem{Tisbury_etal2020b}
 \leavevmode\vrule height 2pt depth -1.6pt width 23pt, {\em The evolution of
   traveling waves in a {K}{P}{P} reaction–diffusion model with cut-off
   reaction rate. {I}{I}. {E}volution of traveling waves}, Stud. in App. Math.,
   146 (2021), pp.~330--370.

 \bibitem{TzellaVanneste2019}
 {\sc A.~Tzella and J.~Vanneste}, {\em Chemical front propagation in periodic
   flows: {F}{K}{P}{P} versus {G}}, SIAM J. Appl. Math., 79 (2019),
   pp.~131--152.

 \bibitem{VanDyke1975}
 {\sc M.~{Van Dyke}}, {\em Perturbation Methods in Fluid Mechanics}, Parabolic
   Press, 1975.

 \bibitem{Xin2000a}
 {\sc J.~Xin}, {\em Front propagation in heterogeneous media}, SIAM Rev., 42
   (2000), pp.~161--230.

 \bibitem{XinYu2013}
 {\sc J.~Xin and Y.~Yu}, {\em Sharp asymptotic growth laws of turbulent flame
   speeds in cellular flows by inviscid {H}amilton--{J}acobi models}, Ann. I. H.
   Poincar\'e AN, 30 (2013), pp.~1049 -- 1068.

 \end{thebibliography}
 \providecommand{\noopsort}[1]{}\providecommand{\singleletter}[1]{#1}%

\end{document}